\documentclass{article}

\usepackage[english]{babel}
\usepackage{amsmath,amsfonts,amssymb,amsthm, bm}
\usepackage{algorithm}
\usepackage{algpseudocode}
\usepackage{caption}
\usepackage{subcaption}
\usepackage{diagbox, xcolor}

\usepackage[letterpaper,top=2cm,bottom=2cm,left=3cm,right=3cm,marginparwidth=1.75cm]{geometry}

\usepackage{amsmath}
\usepackage{graphicx}
\usepackage[colorlinks=True, allcolors=blue]{hyperref}
\usepackage{authblk}

\title{Graph-Based Convexification of Nested Signal Temporal Logic Constraints for Trajectory Optimization}

\date{}
\author[1, 3]{Thomas Claudet\thanks{Corresponding author: \texttt{thomas.claudet@isae-supaero.fr}}}
\author[1]{Davide Martire}
\author[2, 3]{Damiana Losa} 
\author[1]{Francesco Sanfedino}
\author[1]{Daniel Alazard}
\affil[1]{\small{ISAE-SUPAERO, 10 Av. E. Belin, 31400 Toulouse, France}}
\affil[2]{Thales Alenia Space, 5 All. des Gabians, 06150 Cannes, France}
\affil[3]{IRT Saint-Exupéry, 3 Rue Tarfaya, 31400 Toulouse, France}

\begin{document}
\maketitle

\begin{abstract}
Optimizing high-level mission planning constraints is traditionally solved in exponential time and requires to split the problem into several ones, making the connections between them a convoluted task. This paper aims at generalizing recent works on the convexification of Signal Temporal Logic (STL) constraints converting them into linear approximations. Graphs are employed to build general linguistic semantics based on key words (such as \textit{Not}, \textit{And}, \textit{Or}, \textit{Eventually}, \textit{Always}), and super-operators (e.g., \textit{Until}, \textit{Implies}, \textit{If and Only If}) based on already defined ones. Numerical validations demonstrate the performance of the proposed approach on two practical use-cases of satellite optimal guidance using a modified Successive Convexification scheme.

\end{abstract}

\section{Introduction}

The widespread approach for optimizing high-level mission planning scenarios generally involves parallel or distributed frameworks requiring to split the problem into several segments, often failing in finding global optimal solutions \cite{yang2019}. Tasks requiring fast generations may opt for fewer segments and more lenient splitting tolerances in their planning, hence trading for precision. On the other hand, highly accurate trajectories involve increased subproblems with tighter splitting tolerances, leading to longer planning times \cite{wang2021}. In any case, tuning of hyper-parameters (number of segments, tolerance) is an inherent limitation which has to be differently approached from one problem to another. By treating the trajectory as a whole, this issue is alleviated, but at the cost of even larger computational efforts. Generally solved using Mixed-Integer Linear Programming (MILP) or related NP-hard techniques \cite{richards2005}, these methods suffer from poor scaling capabilities and their integration into agile autonomous systems with low computational capabilities today still remains out of reach. \\

Another way of looking at high-level mission planning is from the perspective of logic operators. Derived from the linear temporal logic (LTL) theory \cite{pnueli1977}, logical operators such as \textit{And, Or, Eventually, Always, Until} provide a formal language to express complex temporal constraints and objectives in optimization models. In \cite{donze2010}, Signal Temporal Logic (STL) was introduced as a derivative of LTL to guide simulation-based verification of complex nonlinear and hybrid systems against temporal properties. STL was later extended from verification to optimization to fit within a robust Model Predictive Control (MPC) scheme \cite{Raman2014}. Although enabling more intuitive problem-solving approaches, its underlying MILP baseline showed incompressible limits to its applicability for fast, responsive autonomous systems. \\

To remedy this, in 2022, Y. Mao, B. Açıkmeşe, et al. \cite{scvx_stl} were the first to apply STL to convex optimization \cite{boyd2004}, specifically to the Successive Convexification (SCvx) scheme introduced in \cite{scvx}, and benefit from all its fast and global convergence properties to a local optimum provided convexity assumptions. Signal Temporal Logic was also infused into Machine Learning and Deep Learning \cite{Ma2020}, \cite{vazquez2017} as well as reinforcement learning \cite{Puranic2021}. Furthermore, \cite{leung2023backpropagation} computed the quantitative semantics of STL formulas using computation graphs and perform backpropagation for neural networks using off-the-shelf automatic differentiation tools. However, supervised and reinforcement learning, and especially applied to trajectory optimization, are black boxes requiring enormous amount of realistic data, which makes them problem dependant, training intensive and has no certification for safety critical missions. \\

For the present study, inherited from the foundations laid down in \cite{scvx_stl} and \cite{leung2023backpropagation}, the aim is to bridge the gap between graph-based Signal Temporal Logic and convexification for optimization problem solving. By accompanying the reader through several examples, arbitrarily complex graphs of STL constraints are being discretized and linearized to fit the convex framework. New super-structure of operators are introduced while numerical simulations show the effectiveness of the proposed approach within a modified Successive Convexification scheme. The layout of this work goes as follows:

\begin{itemize}
    \item Section 2 presents the STL grammar and its convexification basics;
\item Section 3 focuses on each basic logical operator individually;
\item Section 4 presents operator nesting and proposes a convexification process based on graphs;
\item Section 5 presents new operators which can be derived from existing ones, therefore creating super-operators and getting closer to real-world operations;
\item Section 6 presents simulation results of optimization problems containing real-world operating constraints;
\item Section 7 finally draws conclusions on this work.

\end{itemize}

\section{STL Grammar and Convexification}\label{sec:prelim}

This section presents the STL grammar (see also \cite{scvx_stl}) and the associated convexification process in the context of a general convex optimization problem of type:
\begin{equation}
    \underset{\boldsymbol{u}}{\mathrm{minimize}} \sum_{k=1}^N J \left(\boldsymbol{x}_k,\boldsymbol{u}_k\right),
\end{equation}
subject to the discrete dynamics:
\begin{equation}
\begin{matrix}
     \boldsymbol{x}_{k+1} = \boldsymbol{f}
    \left(\boldsymbol{x}_k,\boldsymbol{u}_k\right),& k = 1,2,\dotso,N-1,  
\end{matrix}
\end{equation}
as well as STL type constraints being introduced in the present work.
Additional convex constraints of inequality and equality types can also be implemented in the general vector format:
\begin{equation}
    \left\{\begin{array}{ll}
            \boldsymbol{g}(\boldsymbol{x}_k,\boldsymbol{u}_k) \leq \boldsymbol{0}\\
            \boldsymbol{h}(\boldsymbol{x}_k,\boldsymbol{u}_k) = \boldsymbol{0}
    \end{array}\right., k=1,2,\dotso,N.
\end{equation}


\subsection{STL Grammar}

Signal Temporal Logic uses a specific grammar that is being recalled in this section. The following notations are functional to the STL grammar:
\begin{itemize}
    \item $\boldsymbol{x}_k$ and $\boldsymbol{u}_k$: state and control vectors at time instance $k$;
\item $\varphi$ or $\psi$: examples of STL predicates found throughout the paper (always expressed in terms of inequality, an equality being a double inequality);
\item $\mu$: generic function of the state and control vectors;
\item $\pi^\mu$: Boolean predicate whose truth value is determined by the sign of the evaluation of  $\mu$ (i.e. $\pi^\mu =1 \iff \mu > 0$, 0 otherwise).
\end{itemize}

Examples of STL formula could be a function expressing the command variable greater than a certain value or expressing the distance between two objects lower than a threshold value.
Table \ref{table:logics1} presents mathematical key words used in the sequel of this paper and STL operators as building blocks of logical constraints.


\begin{table}[h]

\caption{Mathematical key words and STL operators.}
\quad \quad\quad \quad\quad \quad
\begin{tabular}{ |c|c|} 
\hline
\textbf{Mathematical Symbol} & \textbf{Meaning}\\
 \hline
 Negation & $\neg$ \\
 \hline
 There exists & $\exists$ \\
 \hline
 For all & $\forall$ \\
 \hline
 Implies & $\implies$ \\ 
 \hline
 If and Only If & $\iff$ \\ 
 \hline
Satisfies & $\models$ \\ 
 \hline
 Is Equivalent & $\equiv$ \\ 
 \hline
\end{tabular}
\label{table:logics1}
\quad \quad
\quad \quad
\begin{tabular}{ |c|c|} 
\hline
\textbf{STL Operator} & \textbf{Symbol} \\
 \hline
 \textit{And} & $\wedge$ \\
 \hline
 \textit{Or} & $\lor$ \\
 \hline
 \textit{Eventually} & $\diamond$ \\
 \hline
 \textit{Always} & $\Box$ \\
 \hline
  \textit{Exclusive Or} & $\oplus$ \\
 \hline
 \textit{Until} & $\mathcal{U}$ \\ 
 \hline
\end{tabular}
\label{table:logics2}
\end{table}

The formal semantics of predicates $\varphi$ or/and $\psi$ with respect to a trajectory $\xi$ are defined as follows:
\begin{itemize}
 
 \item A trajectory $\xi$ at time $t_k$, denoted $(\xi, t_k)$ satisfies a predicate $\pi^\mu$ if and only if the associated function $\mu$ at time $k$ is positive:

   \begin{equation}\label{eq: STL def mu}
            (\xi, t_k) \models \pi^\mu \iff \mu(\boldsymbol{x}_k, \boldsymbol{u}_k) > 0.
        \end{equation}

 \item A trajectory $\xi$ at time $t_k$ satisfies the negation of a predicate $\varphi$ if and only if the opposite trajectory satisfies $\varphi$:

   \begin{equation}\label{eq: STL def neg}
             (\xi, t_k) \models \neg \varphi \iff \neg  (\xi, t_k) \models \varphi.
        \end{equation}
 
\item A trajectory $\xi$ at time $t_k$ satisfies the conjunction $\varphi \wedge \psi$ 
 (respectively disjunction $\varphi \lor \psi$) if and only if the signal $\xi$ at time $t_k$ satisfies the STL predicates $\varphi$ and $\psi$ (respectively $\varphi$ or $\psi$):

 \begin{equation}\label{eq: STL def and}
            (\xi, t_k) \models \varphi \wedge \psi \iff (\xi, t_k) \models \varphi \wedge (\xi, t_k) \models \psi,
        \end{equation}

         \begin{equation}\label{eq: STL def or}
            (\xi, t_k) \models \varphi \lor \psi \iff (\xi, t_k) \models \varphi \lor (\xi, t_k) \models \psi.
        \end{equation}

\item A trajectory $\xi$ satisfies \textit{Eventually} ($\diamond$)  $\varphi$ between time $a$ and $b$ if and only if there exists at least one time in that interval $[a,b]$ where the predicate $\varphi$ is satisfied:

\begin{equation}\label{eq: STL def eventually}
            \xi \models \diamond_{[a, b]} \varphi \iff \exists t_{k} \in [a, b], (\xi, t_{k}) \models \varphi. 
        \end{equation}
 
\item A trajectory $\xi$ satisfies Always ($\Box$) $\varphi$ between time $a$ and $b$ if and only if the predicate $\varphi$ is satisfied for each time between the interval $[a,b]$:  

  \begin{equation}\label{eq: STL def always}
            \xi \models \Box_{[a, b]} \varphi \iff \forall t_{k} \in [a, b], (\xi, t_{k}) \models \varphi. 
        \end{equation} 
        

\end{itemize}

\subsubsection{STL Margin Measure}\label{robust logic}

\begin{figure}[ht]
    \centering
    \includegraphics[width=0.75\linewidth]{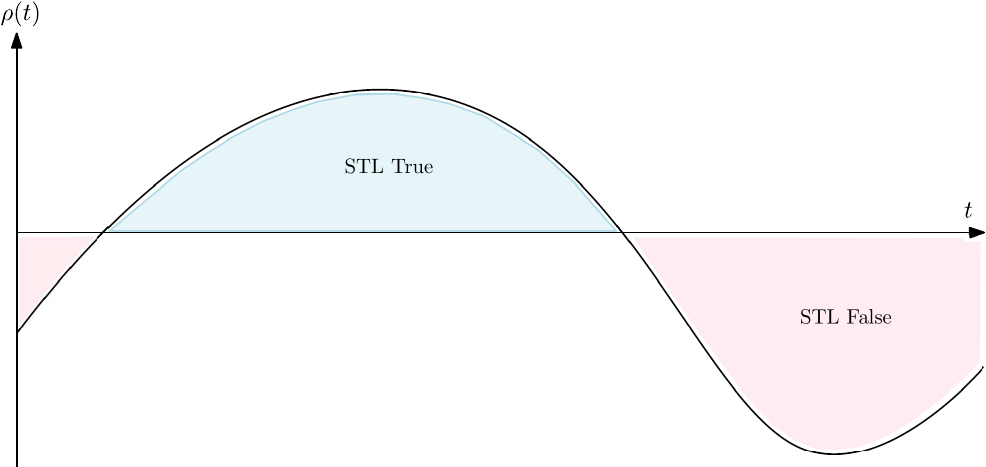}
    
    \caption{The margin function. When the area is blue (respectively pink) the constraint is \textit{True} (respectively \textit{False}).}\label{fig:rho_sign}
\end{figure}

In the above STL Eqs. \eqref{eq: STL def mu} to \eqref{eq: STL def always}, the result is a Boolean: either the property is \textit{True}, or \textit{False}. In order to gain more insight, a margin function is introduced to better classify an STL expression. The margin function $\rho$ defines the truth of the STL formula as shown in Fig. \ref{fig:rho_sign}. This margin function is computed using minimum and maximum functions for the STL operators introduced in Table \ref{table:logics2} and used in Eqs. \eqref{eq: STL def mu} to \eqref{eq: STL def always}. \\

Let us consider a practical example to better illustrate the interest of the margin function $\rho$. By denoting the position of an autonomous system as $\boldsymbol{r}$, the condition $\lVert \boldsymbol{r} (t) \rVert \geq d$ states that at time $t$ the system remains at a distance greater than a certain value $d$ from an obstacle located at the reference frame's origin. In that case, defining a margin function such that $\rho(t) = \lVert \boldsymbol{r}(t) \rVert - d$, the security condition is satisfied at time t if and only if $\rho$ is positive at $t$. STL formulas \eqref{eq: STL def mu} to \eqref{eq: STL def always} can be expressed as follows, in terms of margin semantics using margin functions.\\ 

The complete robust semantics for Eq. \eqref{eq: STL def mu} is given by:

    \begin{equation}
        \rho^{\pi^\mu}_k = \mu(\boldsymbol{x}_k, \boldsymbol{u}_k).
    \end{equation}
Margin function for the negation of a property $\psi$ of the Eq. \eqref{eq: STL def neg} is given by: 
    \begin{equation}
        \rho^{\neg \psi}_k = - \rho^{\psi}_k.
    \end{equation}
     Margin function of Eq. \eqref{eq: STL def and} as conjunction satisfies the following minimum property in terms of margin functions of $\varphi$ and $\psi$: 
    \begin{equation}
        \rho^{\varphi \wedge \psi}_k = \min(\rho^{\varphi}_k, \rho^{\psi}_k).
    \end{equation}
    Instead, margin function of Eq. \eqref{eq: STL def or} as disjunction satisfies the following maximum property in terms or margin functions of $\varphi$ and $\psi$:
    \begin{equation}
        \rho^{\varphi \lor \psi}_k =  \max(\rho^{\varphi}_k, \rho^{\psi}_k).
    \end{equation}
   Margin function of Eq. \eqref{eq: STL def eventually} (involving constraint \textit{Eventually} ($\diamond$) $\varphi$ \textit{True} between time $a$ and $b$) is positive when the maximum of the constraint margin function is positive in the interval: 

    \begin{equation}
        \rho^{\diamond_{[a, b]} \varphi} = \max \limits_{k \in [a, b]}\rho^{\varphi}_k.
    \end{equation}
    Margin function of Eq. \eqref{eq: STL def always} (involving constraint \textit{Always} ($\Box$) $\varphi$ \textit{True} between time $a$ and $b$) is positive when the minimum of the constraint margin function is positive in the interval:
    \begin{equation}
         \rho^ {\square_{[a, b]}{\varphi}}  =  \min \limits_{k \in {[a, b]}}\rho^{\varphi}_k.
    \end{equation}


\subsection{STL Convexification}

The convexification process is made of three steps. First, discretization to make the problem suitable for numerical evaluations and implementation on computers. Then, the second step is linearization to approximate nonlinear functions and work with convenient linear methods. Finally, the third step is to concatenate the linearizations into a linear system of equations, being a suitable input for optimization problem solvers. This section discusses these three steps in greater details.

\subsubsection{Discretization}
As it was seen in Section  \ref{robust logic}, the margin of the operators are functions of time and derived from $min$ and $max$ functions. So, to model this in the discrete world, new variables (scalars), denoted as STL-variables, representing the margin (which is a scalar) evaluated at each time step, are introduced: 
\begin{equation}
    \boldsymbol{\alpha} = \left[\begin{array}{c}
    \alpha_1\\
    \vdots\\
    \alpha_N
    \end{array}\right].
\end{equation}
The STL-variable vector augments the state vector of the optimization variables:
\begin{equation}
    \boldsymbol{x} = \left[\begin{array}{c}
    \boldsymbol{x}_1\\
    \vdots\\
    \boldsymbol{x}_N
    \end{array}\right],
\end{equation}
to give:
\begin{equation}
    \boldsymbol{z} = \left[\begin{array}{c}
    \boldsymbol{x}\\ \hline \boldsymbol{\alpha}
    \end{array}\right].   
\end{equation}

The goal is to replicate the behavior of the \textit{max} or \textit{min} functions over time intervals $[a, b]$. 
The extremum is initialized at the first value of the interval and propagates through it. The obtained value at the end of the interval is kept until the end. It follows:
    \begin{equation}
        \left\{
        \begin{matrix}\label{alpha dyn}
            \begin{array}{ll}
                \alpha_k = \rho^{\varphi}_{k_a}\\
                \alpha_{k+1} = \chi(\alpha_k, \rho^\varphi_{k+1})\\
                \alpha_k = \alpha_{k_b}
            \end{array} &
            \begin{array}{ll}
                \forall k = 1, \dots, k_a\\
                \forall k = k_a, \dots, k_b - 1\\
                \forall k = k_b + 1, \dots, N
            \end{array}
        \end{matrix}
        \right. ,
    \end{equation}
where $\chi$ is a general function representing the margin, e.g., functions $min$ or $max$.
\\
    
Between times $a$ and $b$, the STL-variable vector $\boldsymbol{\alpha}$ is now defined as a sequence of scalar variables $(\alpha_k)_{(k_a:k_b)}$ following the discrete dynamics:
\begin{equation}
\alpha_{k+1} = \chi(\alpha_k, \rho^\varphi_{k+1}).
\end{equation}
The margin function, in turn, depends on the dynamics of the system under control:
\begin{equation}
\rho^\varphi_{k+1}=\rho^\varphi_{k+1}\left(\boldsymbol{x}_{k+1}\right).
\end{equation}

            
\subsubsection{Linearization}\label{subsubsec:linearization}
Now, the STL-variables are defined over the whole time horizon to follow the continuous definition of the margin function. In particular, the functions $min$ and $max$ are not linear as function of the other STL-variables, not even differentiable (indeed the maximum and minimum function have sharp edges at the flips). Function $\rho^\varphi$ is in general not linear as well (e.g. when asking for a distance between an object and an obstacle being greater that a threshold, the norm function is involved). These kinds of functions are not straightforwardly implementable in a standard convex optimization algorithm. A linearization step is then performed to cope with this framework.  \\

Let us assume that function $\chi$ is differentiable. In that case, it would be desired to consider the STL-variables appearing linearly to build a linear system of equations. Using the chain rule (derivation of function composition), the approximation to the first order of the Taylor expansion leads:
\begin{equation}
    \alpha_{k+1} =  \chi(\alpha_k, \rho^\varphi_{k+1}) \approx \chi(\bar{\alpha}_k, \bar{\rho}^\varphi_{k+1}) + \left. \frac{\partial \chi}{\partial \alpha_k}  \right |_{\bar{\alpha}_k, \bar{\rho}^\varphi_{k+1}} (\alpha_k - \bar{\alpha}_k) + \left.\frac{\partial \chi}{\partial \rho^\varphi_{k+1}}  \right |_{\bar{\alpha}_k, \bar{\rho}^\varphi_{k+1}} \left. \frac{\partial \rho^\varphi_{k+1}}{\partial \boldsymbol{x}_{k+1}} \right |_{\bar{\boldsymbol{x}}_{k+1}} (\boldsymbol{x}_{k+1} - \bar{\boldsymbol{x}}_{k+1}),
\end{equation} 
where the notation $\bar{x}_j$ for any general variable $x$ represents the reference of $x$ at time $j$ (e.g. either of the initialization or the last optimal iteration). Introducing the  zero order term: 
\begin{equation}
    R(\bar{\alpha}_k, \bar{\boldsymbol{x}}_{k+1}) = \chi(\bar{\alpha}_k, \bar{\rho}^\varphi_{k+1}) -  \left. \frac{\partial \chi}{\partial \alpha_k}  \right |_{\bar{\alpha}_k, \bar{\rho}^\varphi_{k+1}} \bar{\alpha}_k - \left.\frac{\partial \chi}{\partial \rho^\varphi_{k+1}}  \right |_{\bar{\alpha}_k, \bar{\rho}^\varphi_{k+1}} \left. \frac{\partial \rho^\varphi_{k+1}}{\partial \boldsymbol{x}_{k+1}} \right |_{\bar{\boldsymbol{x}}_{k+1}} \bar{\boldsymbol{x}}_{k+1},
\end{equation} 
Dropping the reference notation for the gradients in the following part of the paper, for the sake of conciseness, one finally gets:
\begin{equation} \label{zeta def}
    \alpha_{k+1} =  R(\bar{\alpha}_k, \bar{\boldsymbol{x}}_{k+1}) + \frac{\partial \chi}{\partial \alpha_k} \alpha_k + \frac{\partial \chi}{\partial \rho^\varphi_{k+1}} \frac{\partial \rho^\varphi_{k+1}}{\partial \boldsymbol{x}_{k+1}}\boldsymbol{x}_{k+1}.
\end{equation}     
        
\begin{figure}[ht]
        \centering
        \includegraphics[width=0.6\linewidth]{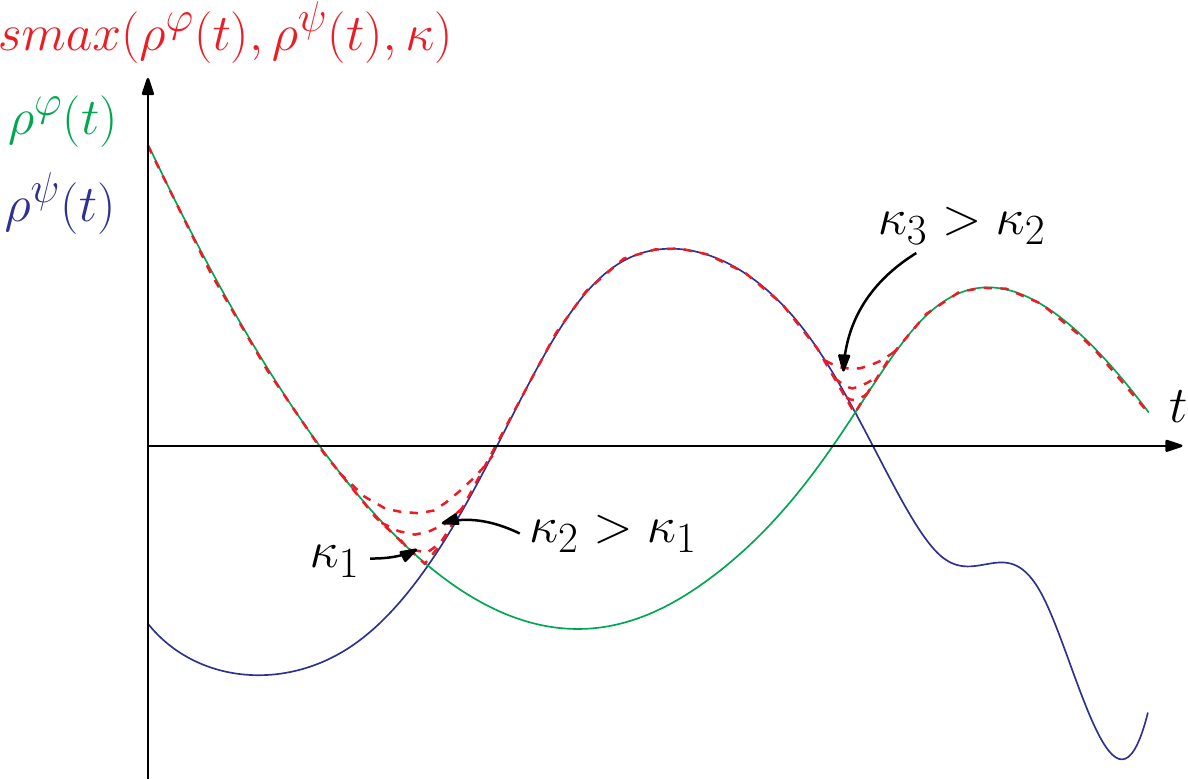}
        
        \caption{The smoothed margin function with respect to time. The green (respectively blue) line represents the evolution of $\rho^\varphi$ (respectively $\rho^\psi$). The red dashed lines represent different smoothing parameters when using the $smax$ function.}\label{fig:smax}
    \end{figure}

Equation \eqref{zeta def} is now linear in the STL-variables $\boldsymbol{\alpha}$ and in the dynamic state variables $\boldsymbol{z}$. Function $\chi$, when equal to $min$ or $max$, is not differentiable. However, as proposed in \cite{scvx_stl}, it can be approximated with smoothed functions ($smin$ and $smax$ for $min$ and $max$ respectively). Function $smin$ and $smax$ are defined as follows (see also Fig. \ref{fig:smax} depicting an example):
        \begin{equation}
            smin(a, b, \kappa) = \left\{
            \begin{matrix}
                \begin{array}{ll}
                        b\\a\\g_{min}(a, b, \kappa)
                \end{array} &
                \begin{array}{ll}
                        \mathrm{if}\\\mathrm{if}\\\mathrm{if}
                \end{array} &
                \begin{array}{ll}
                        a - b \geq \kappa\\
                        a - b \leq -\kappa\\
                        a - b \in [-\kappa, \dots, \kappa]
                \end{array}                 
            \end{matrix}
            \right. ,
        \end{equation}        
        where $g_{min}$ is the smoothing function:
        \begin{equation}
           \begin{matrix}
            g_{min}(a, b, \kappa) = a(1 - h) + hb - \kappa h(1-h), &\mathrm{with} & h = \frac{1}{2} + \frac{a-b}{2\kappa}
            \end{matrix}, 
        \end{equation}
        dependent on its smoothing gain $\kappa$.\\

        As stated in Eq. \eqref{zeta def}, partial derivatives of the margin functions are essential to define the STL-variable dynamics:
        \begin{equation}
            \frac{\partial smin}{\partial a}(a, b, \kappa) = \left\{
            \begin{matrix}
                \begin{array}{ll}
                    0\\1\\1-h
                \end{array} &
                \begin{array}{ll}
                    \mathrm{if}\\\mathrm{if}\\\mathrm{if}
                \end{array} &
                \begin{array}{ll}
                    a - b \geq \kappa\\
                    a - b \leq -\kappa\\
                    a - b \in [-\kappa, \dots, \kappa]
                \end{array} 
            \end{matrix}
            \right.,
        \end{equation}
        \begin{equation}
            \frac{\partial smin}{\partial b}(a, b, \kappa) = \left\{
            \begin{matrix}
                \begin{array}{ll}
                    1\\0\\h
                \end{array} &
                \begin{array}{ll}
                    \mathrm{if}\\\mathrm{if}\\\mathrm{if}
                \end{array} &
                \begin{array}{ll}
                    a - b \geq \kappa\\
                    a - b \leq -\kappa\\
                    a - b \in [-\kappa, \dots, \kappa]
                \end{array} 
            \end{matrix} 
            \right. .
        \end{equation}

        Caution has to be taken by the reader since $smin \neq -smax$. Indeed, function $smax$ is defined as:

        \begin{equation}
            smax(a, b, \kappa) = \left\{
            \begin{matrix}
                \begin{array}{ll}
                        b\\a\\g_{max}(a, b, \kappa)
                \end{array} &
                \begin{array}{ll}
                        \mathrm{if}\\\mathrm{if}\\\mathrm{if}
                \end{array} &
                \begin{array}{ll}
                        b - a \geq \kappa\\
                        b - a \leq -\kappa\\
                        b - a \in [-\kappa, \dots, \kappa]
                \end{array}                 
            \end{matrix}
            \right. ,
        \end{equation}     
        where $g_{max}$ is the smoothing function:
         \begin{equation}
           \begin{matrix}
            g_{max}(a, b, \kappa) = b(1 - h) + ha + \kappa h(1-h), &\mathrm{with} & h = \frac{1}{2} + \frac{a - b}{2\kappa}
            \end{matrix},
        \end{equation}
        dependent on its smoothing gain $\kappa$.\\
        
        Partial derivatives of the margin function $smax$ are the following:
        \begin{equation}
            \frac{\partial smin}{\partial a}(a, b, \kappa) = \left\{
            \begin{matrix}
                \begin{array}{ll}
                    0\\1\\h
                \end{array} &
                \begin{array}{ll}
                    \mathrm{if}\\\mathrm{if}\\\mathrm{if}
                \end{array} &
                \begin{array}{ll}
                    b - a \geq \kappa\\
                    b - a \leq -\kappa\\
                    b - a \in [-\kappa, \dots, \kappa]
                \end{array} 
            \end{matrix}
            \right.,
        \end{equation}
        \begin{equation}
            \frac{\partial smin}{\partial b}(a, b, \kappa) = \left\{
            \begin{matrix}
                \begin{array}{ll}
                    1\\0\\1-h
                \end{array} &
                \begin{array}{ll}
                    \mathrm{if}\\\mathrm{if}\\\mathrm{if}
                \end{array} &
                \begin{array}{ll}
                    b - a \geq \kappa\\
                    b - a \leq -\kappa\\
                    b - a \in [-\kappa, \dots, \kappa]
                \end{array} 
            \end{matrix} 
            \right. .
        \end{equation}

         These partial derivatives will allow to implement the derivatives of $\chi = \{min, max\}$ with respect to its variables as seen in \eqref{zeta def}.\\
         
        In cases where the margin functions might not be differentiable at certain points (e.g., the norm function is not differentiable at 0), a null gradient can, for instance, be automatically introduced when performing the computation at these singular points. Of course, the smaller the $\kappa$, the better the approximation (in practice, not thinking about differentiability, $\kappa = 0$ seemed of no concern in this work). Maybe some edge cases could be found in highly nonlinear settings, but with a usual step size, the vicinity of the flips is rarely, not to say never, encountered in practical applications. The takeaway is probably that that there is a parameter which can smoothen the dynamics if needed. In the following, the $\kappa$ input of $\chi$ will be omitted for conciseness.

\subsubsection{Convex Framework}\label{cvx framework}

After discretizing and linearizing the robust semantics of the STL-formulas, the final step is to transform the STL-variable dynamics equations into a convex matrix form of type $\boldsymbol{Az} = \boldsymbol{b}$, where:
\begin{itemize}
    \item $\boldsymbol{A}$ is defined as the state-transition matrix;
    \item $\boldsymbol{z}$ is the augmented optimization vector composed of the state optimization variable vector $\boldsymbol{x}$ and of the STL-variable vector $\boldsymbol{\alpha}$);
    \item $\boldsymbol{b}$ is the vector independent on the optimization variables.
\end{itemize}

To define these three convex framework elements on the basis of the discrete and linearized dynamics of Eqs \eqref{alpha dyn} and \eqref{zeta def}, the following mathematical notation are used:
\begin{itemize}
    \item $\boldsymbol{I}_K$ is an identity square matrix of size $K$;
    \item $\boldsymbol{0}_K$ is a column vector of length $K$ with all its elements equal to 0;
    \item $\boldsymbol{0}_{K_1,K_2}$ is a zero matrix with $K_1$ rows and $K_2$ columns;
    \item $\boldsymbol{1}_K$ is a column vector of length $K$ with all its elements equal to 1;
    \item $\mathrm{diag}\left(\left[\begin{matrix} u_1\dots u_m \end{matrix}\right]\right)_K$ is a diagonal concatenation of $K$ blocks, each equal to the row vector $\left[\begin{matrix} u_1\dots u_m\end{matrix}\right]$;
    \item $\mathrm{diag}\left(\left[f(v_{1k}), \dots f(v_{mk})\right]\right)_{K_1:K_2}$ is a diagonal concatenation of $K_2-K_1+1$ blocks, equal to $\left[f(v_{1k}), \dots f(v_{mk})\right]$ for $k=K_1,\dots,K_2$.
\end{itemize}

Compact convex formulation $\boldsymbol{Az} = \boldsymbol{b}$ is based on the following definitions:
\begin{enumerate}
    \item For matrix $\boldsymbol{A}$:
        \begin{equation} \label{eq:matrixA}
        \boldsymbol{A} = \left[\begin{array}{ccc|ccc}
            \boldsymbol{0}_{k_a, k_a} & \boldsymbol{0}_{k_a, k_b - k_a}& \boldsymbol{0}_{k_a, N - k_b} & \boldsymbol{I}_{k_a} & \boldsymbol{0}_{k_a, k_b - k_a} & \boldsymbol{0}_{k_a, N - k_b} \\    
            \boldsymbol{0}_{k_b - k_a, k_a} &  \boldsymbol{A_{22}} &\boldsymbol{0}_{k_b - k_a, N - k_b}& \boldsymbol{0}_{k_b - k_a, k_a} & \boldsymbol{A_{25}} & \boldsymbol{0}_{k_b - k_a, N - k_b} \\
            \boldsymbol{0}_{N-k_b, k_a} & \boldsymbol{0}_{N-k_b, k_b - k_a} & \boldsymbol{0}_{N-k_b, N - k_b} & \boldsymbol{0}_{N-k_b, k_a} & \boldsymbol{0}_{N-k_b, k_b - k_a} & \boldsymbol{A_{36}}
            \end{array} \right],
        \end{equation}
        with:
        \begin{equation}
            \boldsymbol{A_{22}} = \mathrm{diag}\left(\left[\frac{\partial \chi}{\partial \rho^\varphi_{k+1}} \frac{\partial \rho^\varphi_{k+1}}{\partial \boldsymbol{x}^\varphi_{k+1}}\right]\right)_{k_a:k_b-1},
        \end{equation}
        \begin{equation}
            \boldsymbol{A_{25}} = 
            \mathrm{diag}
            \left(\left[\begin{matrix} \frac{\partial \chi}{\partial \alpha_k}&-1\end{matrix}\right]\right)_{k_a:k_b-1},
        \end{equation}
        \begin{equation}
            \boldsymbol{A_{36}} = 
            \mathrm{diag}
            \left(\left[\begin{matrix} 1&-1\end{matrix}\right]\right)_{N-k_b}.
        \end{equation}
    \item For vector $\boldsymbol{z}$:
        \begin{equation}
        \boldsymbol{z} = \left[\begin{array}{c}
        \boldsymbol{x}^\varphi\\ \hline \boldsymbol{\alpha}
        \end{array}\right].   
        \end{equation}
    \item For vector $\boldsymbol{b}$:
    \begin{equation}\label{eq:vectorb}
        \boldsymbol{b} =\begin{bmatrix} 
        \bar{\rho}^{\varphi}_{k_a}\boldsymbol{1}_{k_a} \\ 
        -\mathrm{diag}\left(\left[R(\bar{\alpha}_k, \bar{\boldsymbol{x}}^\varphi_{k+1})\right]\right)_{k_a: k_b-1}
        \boldsymbol{1}_{k_b-k_a}\\
        \boldsymbol{0}_{N-k_b}
        \end{bmatrix}.
    \end{equation}
\end{enumerate}

This linear system of equations can then be concatenated with the other physical dynamics linear constraints to obtain an augmented optimization problem, which will then be fed to the convex optimization solver (e.g., second-order cone, semi-definite).\\

Finally, it remains a very last constraint, which is actually the most natural and important. Let us recall that the STL constraint will be satisfied if and only if the margin function is positive. This can be enforced by asking the last time step of the STL variable to be positive:
\begin{equation}
    \alpha_N \geq 0.
\end{equation}
Note that all the constraints are made to specifically arrive to this one.
\section{Isolated STL Operators}\label{sec:isolated_operators}
This section presents the STL operators one by one as introduced in \cite{scvx_stl}. In Section \ref{sec:nesting}, nesting of these STL operators considered as building blocks will be presented. This paper introduces the classification of the four main operators \textit{And}, \textit{Or}, \textit{Eventually}, \textit{Always} as operators of type \textit{Bridge Operator} or \textit{Flow}. 
\begin{itemize}
    \item An operator is of \textit{Bridge} type when it connects two inputs and generates only one output. \textit{And} and \textit{Or} are examples of \textit{Bridge} operators.

    \item An operator is of \textit{Flow} type when it is connected to a single input and generates a single output. \textit{Eventually} and \textit{Always} are examples of \textit{Flow} operators.

\end{itemize}

\subsection{\textit{Bridge} Operators}

\subsubsection{\textit{And}}
    
\begin{figure}[ht]
        \centering
        \includegraphics[width=0.2\linewidth]{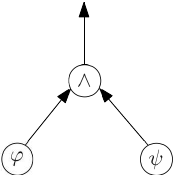}
        \caption{Graph of the \textit{And} operator.}
        \label{fig:and}
    \end{figure}

Let $\varphi$ and $\psi$, be two predicates. The conjunction operator between $\varphi$ and $\psi$ is shown in Fig. \ref{fig:and}, represented by a graph.

In terms of robust STL semantics, the margin function of the conjunction operator can be expressed as a minimum function:
\begin{equation}
    min(\rho^{\varphi}, \rho^{\psi}).
\end{equation}
The discretization gives
\begin{equation}
    \alpha_{k+1} = min(\rho^\varphi_{k+1},\rho^\psi_{k+1}), \;\;\; \forall k = 0, 1, \dots, N-1. 
\end{equation}

Defining as $\boldsymbol{x}^\varphi$ (respectively $\boldsymbol{x}^\psi$) the state variable vector used to write predicate $\varphi$ (respectively $\psi$), the linearization provides:
\begin{equation}
    \alpha_{k+1} = R(\bar{\boldsymbol{x}}^\varphi_{k+1}, \bar{\boldsymbol{x}}^\psi_{k+1}) + \frac{\partial smin}{\partial \rho^{\varphi}_{k+1}}\frac{\partial \rho^{\varphi}_{k+1}}{\partial \boldsymbol{x}^\varphi_{k+1}}\; \boldsymbol{x}^\varphi_{k+1} + \frac{\partial smin}{\partial \rho^{\psi}_{k+1}}\frac{\partial \rho^{\psi}_{k+1}}{\partial \boldsymbol{x}^\psi_{k+1}}\boldsymbol{x}^\psi_{k+1}.
\end{equation}
In general, the states considered by the two margin functions are completely different. For instance, one could use the positions, while the other one the attitude. In that case, they would have to be positioned differently in the STL matrix $\boldsymbol{A}$, easily made with an automated routine.\\

The terms of the compact convex formulation $\boldsymbol{Az} = \boldsymbol{b}$ becomes:


\begin{enumerate}
    \item For matrix $\boldsymbol{A}$:
        \begin{equation}
        \boldsymbol{A} = \begin{bmatrix}\begin{array}{c|c|c}
        \mathrm{diag}\left(\left[\frac{\partial smin}{\partial \rho^{\varphi}_{k+1}}\frac{\partial \rho^{\varphi}_{k+1}}{\partial \boldsymbol{x}^\varphi_{k+1}}\right]\right)_{0:N-1} &  \mathrm{diag}\left(\left[\frac{\partial smin}{\partial \rho^{\psi}_{k+1}}\frac{\partial \rho^{\psi}_{k+1}}{\partial \boldsymbol{x}^\psi_{k+1}}\right]\right)_{0:N-1}  & -\boldsymbol{I}_{N}
        \end{array}
        \end{bmatrix},
        \end{equation}
    \item For vector $\boldsymbol{z}$:
        \begin{equation}\label{vector z}
        \boldsymbol{z} = \left[\begin{array}{c}
         \boldsymbol{x}^\varphi\\ \hline \boldsymbol{x}^\psi\\ \hline
         \boldsymbol{\alpha}
    \end{array}\right],   
        \end{equation}
    \item For vector $\boldsymbol{b}$:
    \begin{equation}\label{vector b}
        \boldsymbol{b} =
        -\mathrm{diag}\left(\left[R(\bar{\boldsymbol{x}}^\varphi_{k+1}, \bar{\boldsymbol{x}}^\psi_{k+1})\right]\right)_{0:N-1}
        \boldsymbol{1}_{N}.
    \end{equation}
\end{enumerate}

\begin{figure}[ht]
        \centering
        \includegraphics[width=0.3\linewidth]{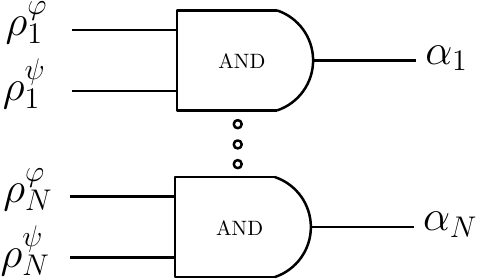}
        \caption{Electronic circuit of the \textit{And} operator.}
        \label{fig:and_elec}
    \end{figure}

Taking the analogy with electric circuits, the system would look like Fig. \ref{fig:and_elec}:\\

\subsubsection{\textit{Or}}   
    
Let $\varphi$ and $\psi$, be two predicates. The disjunction operator between $\varphi$ and $\psi$ is shown in Fig. \ref{fig:and}, represented by a graph.
 \begin{figure}[ht]
        \centering
        \includegraphics[width=0.2\linewidth]{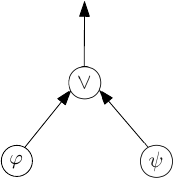}
        \caption{Graph of the \textit{Or} operator.}
        \label{fig:or}
    \end{figure}

In terms of robust STL semantics, the margin function of the disjunction operator can be expressed as a maximum function:
\begin{equation}
    max(\rho^{\varphi}, \rho^{\psi}).
\end{equation}
The discretization gives:
\begin{equation}
    \alpha_{k+1} = max(\rho^\varphi_{k+1},\rho^\psi_{k+1}),\;\;\; \forall k = 0, 1, \dots, N-1.
\end{equation}

Defining as $\boldsymbol{x}^\varphi$ (respectively $\boldsymbol{x}^\psi$) the state variable vector used to write predicate $\varphi$ (respectively $\psi$), the linearization provides:
\begin{equation}
    \alpha_{k+1} = R(\bar{\boldsymbol{x}}^\varphi_{k+1}, \bar{\boldsymbol{x}}^\psi_{k+1}) + \frac{\partial smax}{\partial \rho^{\varphi}_{k+1}}\frac{\partial \rho^{\varphi}_{k+1}}{\partial \boldsymbol{x}^\varphi_{k+1}}\; \boldsymbol{x}^\varphi_{k+1} + \frac{\partial smax}{\partial \rho^{\psi}_{k+1}}\frac{\partial \rho^{\psi}_{k+1}}{\partial \boldsymbol{x}^\psi_{k+1}}\boldsymbol{x}^\psi_{k+1}.
\end{equation}
In general, the states considered by the two margin functions are completely different. For instance, one could use the positions, while the other one the attitude. In that case, they would have to be positioned differently in the STL matrix $\boldsymbol{A}$, easily made with an automated routine.\\

The terms of the compact convex formulation $\boldsymbol{Az} = \boldsymbol{b}$ becomes:

\begin{enumerate}
    \item For matrix $\boldsymbol{A}$:
         \begin{equation}
        \boldsymbol{A} = \begin{bmatrix}\begin{array}{c|c|c}
        \mathrm{diag}\left(\left[\frac{\partial smax}{\partial \rho^{\varphi}_{k+1}}\frac{\partial \rho^{\varphi}_{k+1}}{\partial \boldsymbol{x}^\varphi_{k+1}}\right]\right)_{0:N-1} &  \mathrm{diag}\left(\left[\frac{\partial smax}{\partial \rho^{\psi}_{k+1}}\frac{\partial \rho^{\psi}_{k+1}}{\partial \boldsymbol{x}^\psi_{k+1}}\right]\right)_{0:N-1}  & -\boldsymbol{I}_{N}
        \end{array}
        \end{bmatrix},
        \end{equation}
    \item For vector $\boldsymbol{z}$:
        \begin{equation}
        \boldsymbol{z} = \left[\begin{array}{c}
         \boldsymbol{x}^\varphi\\ \hline \boldsymbol{x}^\psi\\ \hline
         \boldsymbol{\alpha}
    \end{array}\right],   
        \end{equation}
    \item For vector $\boldsymbol{b}$:
    \begin{equation}
        \boldsymbol{b} =
        -\mathrm{diag}\left(\left[R(\bar{\boldsymbol{x}}^\varphi_{k+1}, \bar{\boldsymbol{x}}^\psi_{k+1})\right]\right)_{0:N-1}
        \boldsymbol{1}_{N}.
    \end{equation}
\end{enumerate}

Taking the analogy with electric circuits, the system would look like Fig. \ref{fig:or_elec}:
 \begin{figure}[ht]
        \centering
        \includegraphics[width=0.3\linewidth]{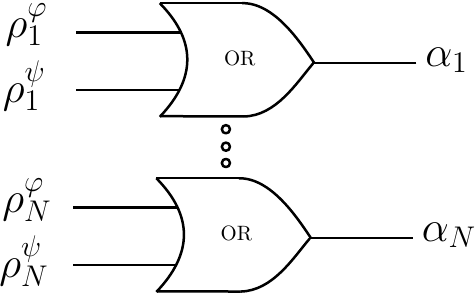}
        \caption{Electronic circuit of the \textit{Or} operator.}
        \label{fig:or_elec}
    \end{figure}

\subsection{\textit{Flow} Operators}

\subsubsection{\textit{Always}}

\begin{figure}[ht]
    \centering
    \includegraphics[width=0.04\linewidth]{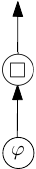}
    \caption{Graph of the \textit{Always} operator.}
    \label{fig:always}
\end{figure}

The \textit{Always} operator is an operator of flow type, it asks for a property to always be \textit{True} inside of a time window, which can cover the whole simulation. The \textit{Always} operator is shown in Fig. \ref{fig:always} represented by a graph. \\

In terms of robust STL semantics, the margin function of this operator can be seen as a recurrent conjunction active at each time step of the time window. 
The discretization of this operator gives:
\begin{equation}
    \left\{
    \begin{matrix}\label{alpha dyn always}
        \begin{array}{ll}
            \alpha_k = \rho^{\varphi}_{k_a}\\
            \alpha_{k+1} = \min(\alpha_k, \rho^\varphi_{k+1})\\
            \alpha_k = \alpha_{k_b}
        \end{array} &
        \begin{array}{ll}
            \forall k = 1, \dots, k_a\\
            \forall k = k_a, \dots, k_b - 1\\
            \forall k = k_b + 1, \dots, N
        \end{array}
    \end{matrix}
    \right. ,
\end{equation}

\begin{figure}[ht]
    \centering
    \includegraphics[width=0.5\linewidth]{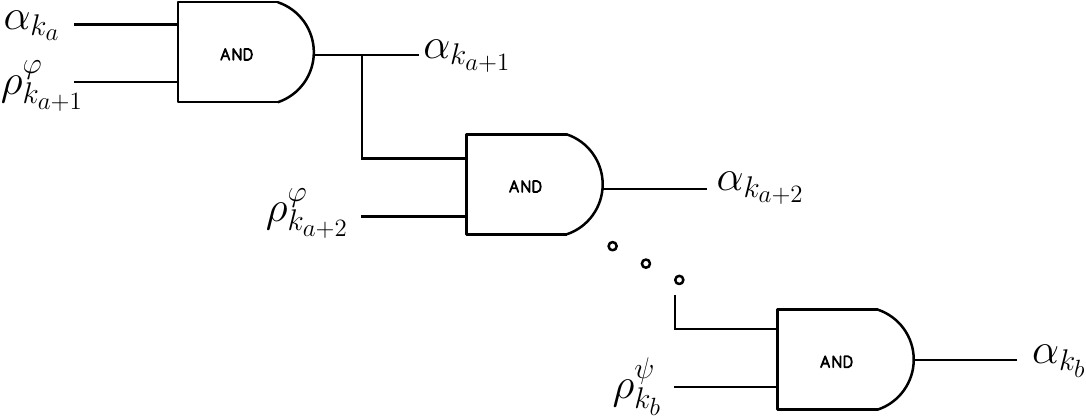}
            \caption{Electronic circuitry of the \textit{Always} operator.}
    \label{fig:always_elec}
\end{figure}

The terms of the compact convex formulation $\boldsymbol{Az} = \boldsymbol{b}$ have the same expression as the ones in Eqs \eqref{eq:matrixA}-\eqref{eq:vectorb} with:

    \begin{equation}
        \boldsymbol{A_{22}} = \mathrm{diag}\left(\left[\frac{\partial smin}{\partial \rho^\varphi_{k+1}} \frac{\partial \rho^\varphi_{k+1}}{\partial \boldsymbol{x}^\varphi_{k+1}}\right]\right)_{k_a: k_b-1},
    \end{equation}
    \begin{equation}
        \boldsymbol{A_{25}} = 
        \mathrm{diag}
        \left(\left[\begin{matrix} \frac{\partial smin}{\partial \alpha_k}&-1\end{matrix}\right]\right)_{k_a:k_b-1}.
    \end{equation}
\\
Taking the analogy with electric circuits, the system would look like Fig. \ref{fig:always_elec}.

\subsubsection{\textit{Eventually}}
\begin{figure}[ht]
        \centering
        \includegraphics[width=0.04\linewidth]{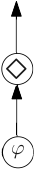}
        \caption{Graph of the \textit{Eventually} operator.}
        \label{fig:eventually}
\end{figure}

The \textit{Eventually} operator is an operator of flow type, it asks for a property to be \textit{True} at least once inside of a time window, which can cover the whole simulation. The \textit{Eventually} operator is shown in Fig. \ref{fig:eventually} represented by a graph. 

In terms of robust STL semantics, the margin function of this operator can be seen as a recurrent disjunction active at each time step of the time window. 
The discretization of this operator gives:
\begin{equation}
    \left\{
    \begin{matrix}\label{alpha dyn eventually}
        \begin{array}{ll}
            \alpha_k = \rho^{\varphi}_{k_a}\\
            \alpha_{k+1} = \max(\alpha_k, \rho^\varphi_{k+1})\\
            \alpha_k = \alpha_{k_b}
        \end{array} &
        \begin{array}{ll}
            \forall k = 1, \dots, k_a\\
            \forall k = k_a, \dots, k_b - 1\\
            \forall k = k_b + 1, \dots, N
        \end{array}
    \end{matrix}
    \right. .
\end{equation}

\begin{figure}[ht]
        \centering
        \includegraphics[width=0.5\linewidth]{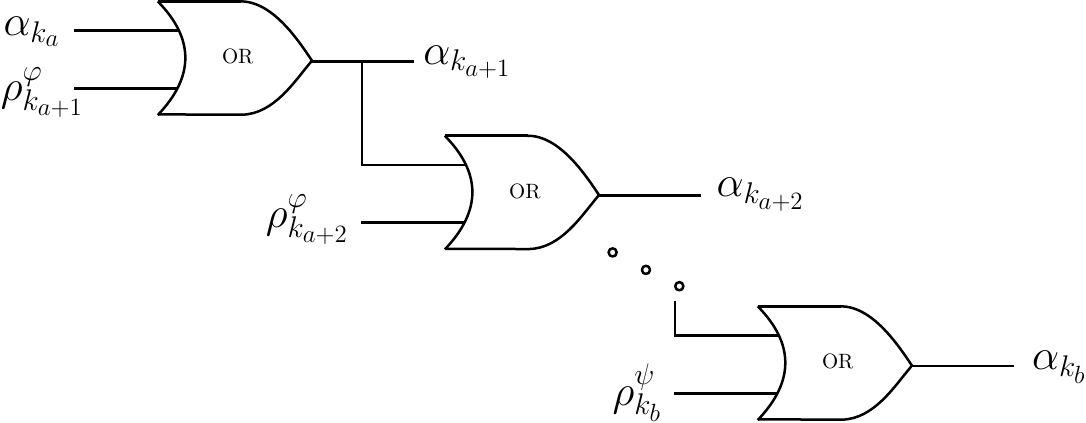}
                \caption{Electronic circuitry of the \textit{Eventually} operator.}
        \label{fig:eventually_elec}
    \end{figure}

The terms of the compact convex formulation $\boldsymbol{Az} = \boldsymbol{b}$ have the same expression of Eqs \eqref{eq:matrixA}-\eqref{eq:vectorb} with:

    \begin{equation}
        \boldsymbol{A_{22}} = \mathrm{diag}\left(\left[\frac{\partial smax}{\partial \rho^\varphi_{k+1}} \frac{\partial \rho^\varphi_{k+1}}{\partial \boldsymbol{x}^\varphi_{k+1}}\right]\right)_{k_a:k_b-1},
    \end{equation}
    \begin{equation}
        \boldsymbol{A_{25}} = 
        \mathrm{diag}
        \left(\left[\begin{matrix} \frac{\partial smax}{\partial \alpha_k}&-1\end{matrix}\right]\right)_{k_a:k_b-1}.
    \end{equation}
\\
Taking the analogy with electric circuits, the system would look like Fig. \ref{fig:eventually_elec}.

\section{Nested STL Operators}\label{sec:nesting}
Nested (i.e., concatenated) STL operators can be used to express complex constraints. This section presents a method to deal with nested STL operators based on graphs. In a graph representing a complex logic expression, the roots (i.e., the first generation parents) are all the predicates of the expression.
The top node of the graph (i.e., the youngest child) must be an operator of \textit{Flow} type (e.g., \textit{Always} or \textit{Eventually}).
This means that, for example, at the very top of the graph, a \textit{Bridge} expression as $\varphi \lor \psi$ is not accepted. To be more precise, a \textit{Bridge} operator always has to be attached to a \textit{Flow} operator above it (e.g., to form $\Box (\varphi \lor \psi))$. This allows to obtain a unequivocal output that does not lead to misunderstandings.
\\

Following this logic, it is possible for each operator to only look at one child down the tree and there is a general formulation to connect them to the rest of the tree.\\

In the following, the naming conventions of STL variables follows the order of the Greek alphabet the higher in the tree the operator is. 



\subsection{Nested \textit{Flow} Operators}

\begin{figure}[H]
    \centering
    \begin{subfigure}[b]{0.25\textwidth}
        \centering
        \includegraphics[width=1.\linewidth]{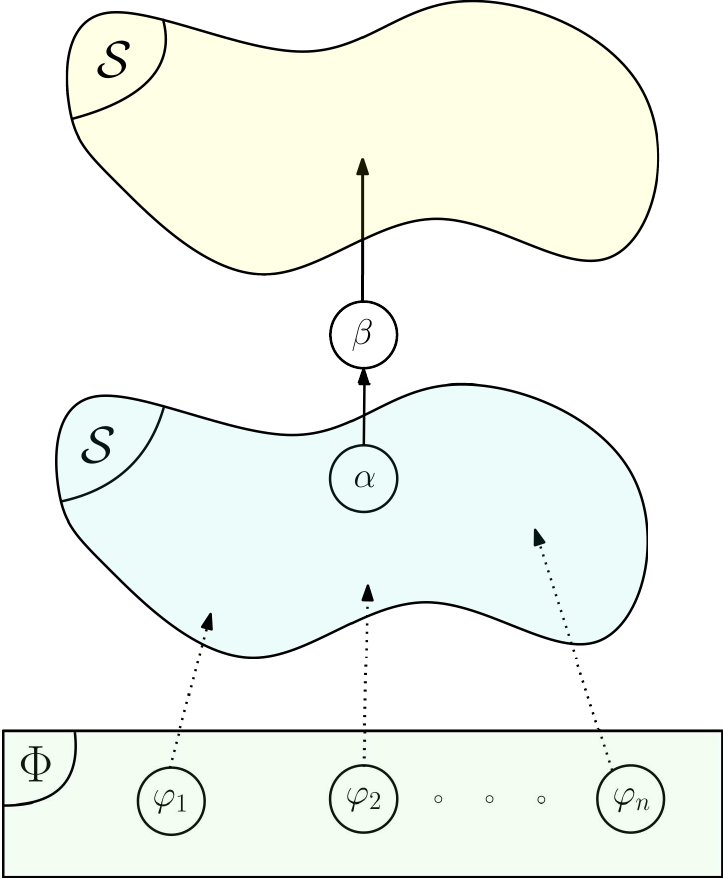}
       \caption{}
        \label{fig:flow1}
    \end{subfigure}
    \;\;\;\;\;\;\;\;\;\;\;\;\;\;\;\;\;\;
    \begin{subfigure}[b]{0.25\textwidth}
        \includegraphics[width=1.\linewidth]{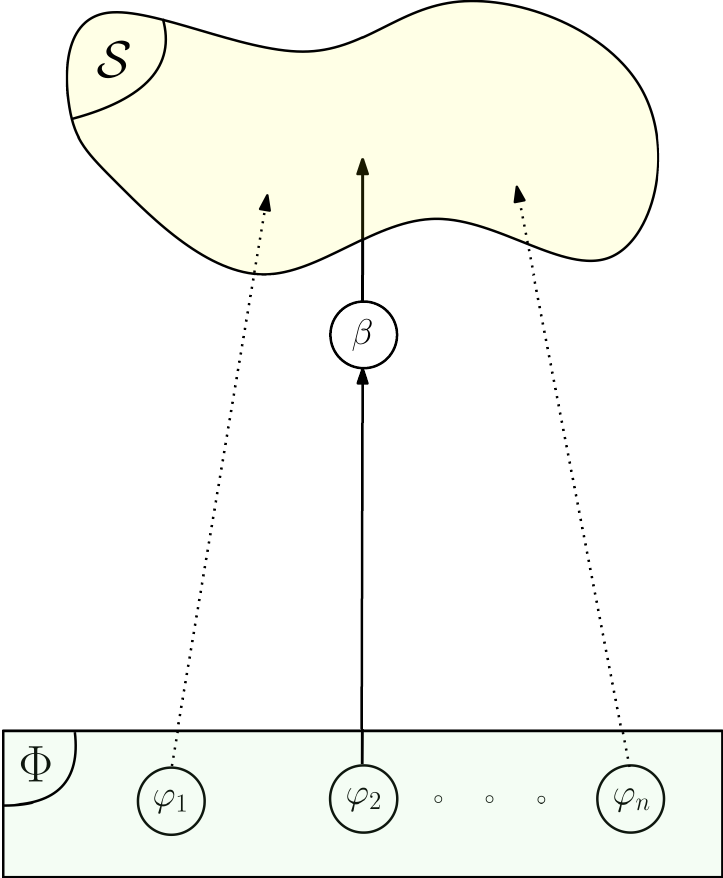}
        \caption{}
        \label{fig:flow2}
    \end{subfigure}
    \caption{General graphs for the \textit{Flow} operators. (a) Connection to another STL variable. (b) Direct connection to a predicate. The clouds in blue and yellow represent the rest of the graph. $\mathcal{S}$ is the STL variable space while $\Phi$ represents the root predicate space.}\label{fig:flow}
\end{figure}

A \textit{Flow} type operator is, by definition, connected to a single input. This is the case for the operator with STL variable vector $\boldsymbol{\beta}$ in Fig. \ref{fig:flow}. The input \textit{Flow} type operator can either be the output of another STL operator (e.g., $\Box$ or $\diamond$), or a simple predicate (e.g., $\varphi$). \\

The equations governing the transition at the \textit{Flow} operator level are the following for the case represented in Fig. \ref{fig:flow1} and \ref{fig:flow2} respectively:
\begin{equation}
    \begin{matrix}
    \beta_{k+1} =
    \frac{\partial \chi}{\partial \beta_k}\beta_k + \frac{\partial \chi}{\partial \alpha_{k+1}}\alpha_{k+1} + R(\bar{\beta}_k, \bar{\alpha} _{k+1}) & \text{if} & \mathcal{P^\beta} \in \mathcal{S},
    \end{matrix}
\end{equation}
\begin{equation}
    \begin{matrix}
    \beta_{k+1} =
    \frac{\partial \chi}{\partial \beta_k}\beta_k + \frac{\partial \chi}{\partial \rho^{\varphi}_{k+1}}\frac{\partial \rho^{\varphi}_{k+1}}{\partial \boldsymbol{x}_{k+1}}\boldsymbol{x}_{k+1} + R(\bar{\beta}_k, \bar{\boldsymbol{x}} _{k+1})& \text{if} & \mathcal{P^\beta} \in \Phi.
    \end{matrix}
\end{equation}

In the above expressions  
\begin{itemize}
    \item $\mathcal{P}^\beta$ is the space of the parent node associated to STL variable $\boldsymbol{\beta}$;
    \item $\mathcal{S}$ is the space of all the STL variables of all the STL complex expressions;
    \item $\Phi$ is the root predicate space.
\end{itemize}
For instance, if the \textit{Flow} operator is connected to an STL expression (e.g., $\Box (\varphi \wedge \psi)$) as in Fig. \ref{fig:flow1}, then $\mathcal{P}^\beta \in \mathcal{S}$. On the other hand, if the flow is directly linked to the predicate $\varphi$ as in Fig. \ref{fig:flow2}, then $\mathcal{P}^\beta \in \Phi$. 
In cases where the time windows do not cover the whole simulations, the three cases as shown in Eq. \eqref{alpha dyn always} or \eqref{alpha dyn eventually} are being used. They are not written again for conciseness. \\


    


\subsection{Nested \textit{Bridge} Operators}

\begin{figure}[H]
    \centering
    \begin{subfigure}[b]{0.24\textwidth}
        \centering
        \includegraphics[width=1\linewidth]{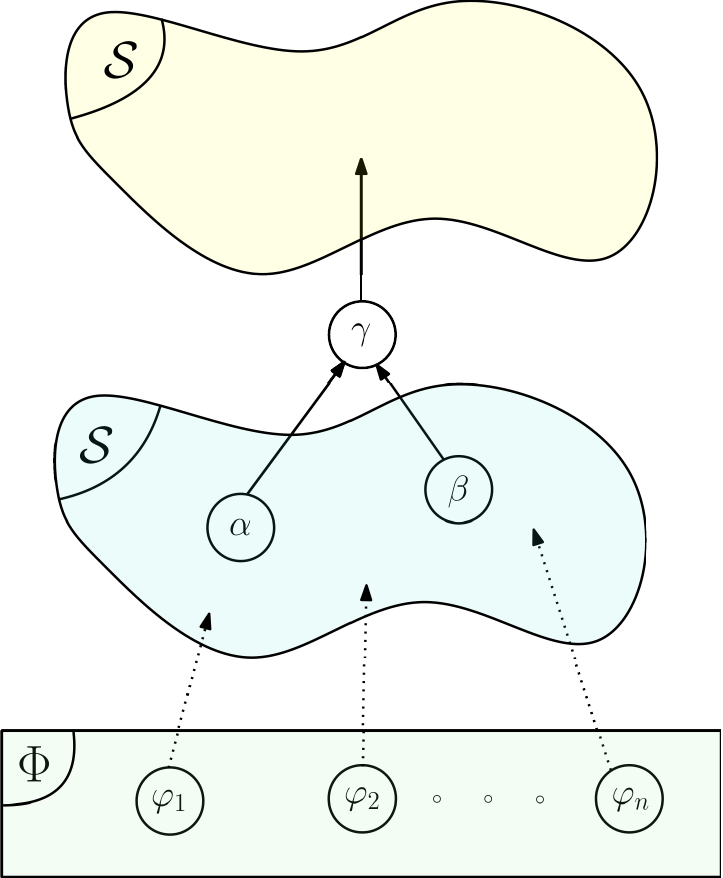}
       \caption{}
        \label{fig:bridge1}
    \end{subfigure}
    \hfill
    \begin{subfigure}[b]{0.24\textwidth}
        \includegraphics[width=1\linewidth]{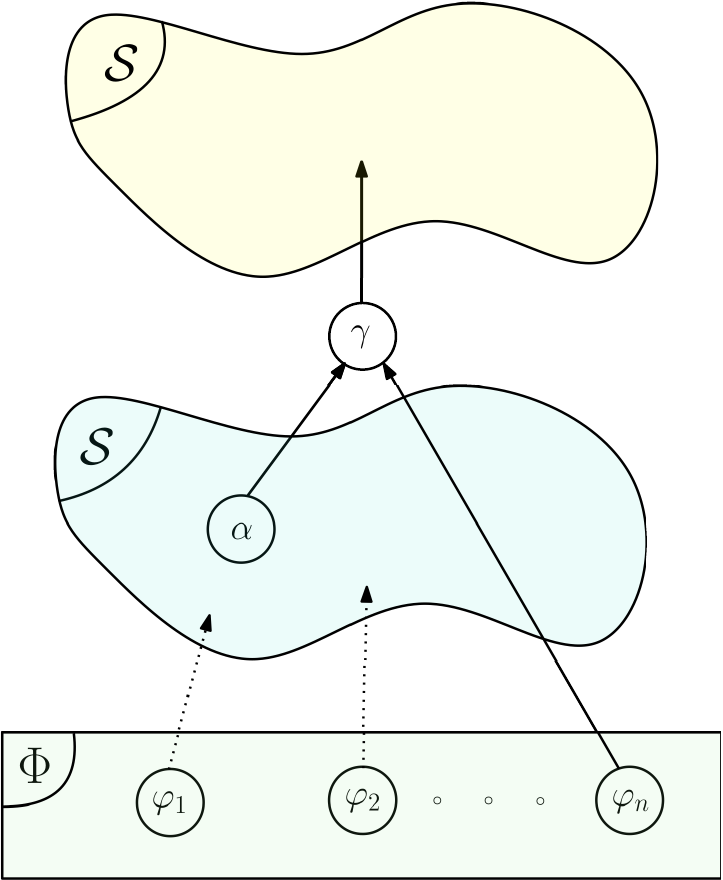}
        \caption{}
        \label{fig:bridge2}
    \end{subfigure}
    \hfill
    \begin{subfigure}[b]{0.24\textwidth}
        \centering
        \includegraphics[width=1\linewidth]{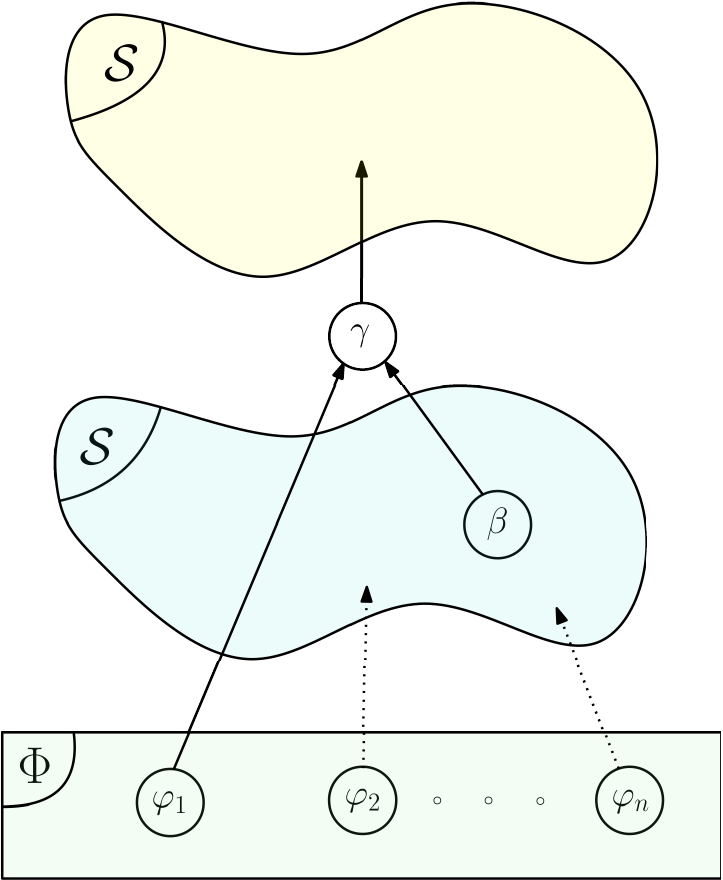}
       \caption{}
        \label{fig:bridge3}
    \end{subfigure}
    \hfill
    \begin{subfigure}[b]{0.24\textwidth}
        \centering
        \includegraphics[width=1\linewidth]{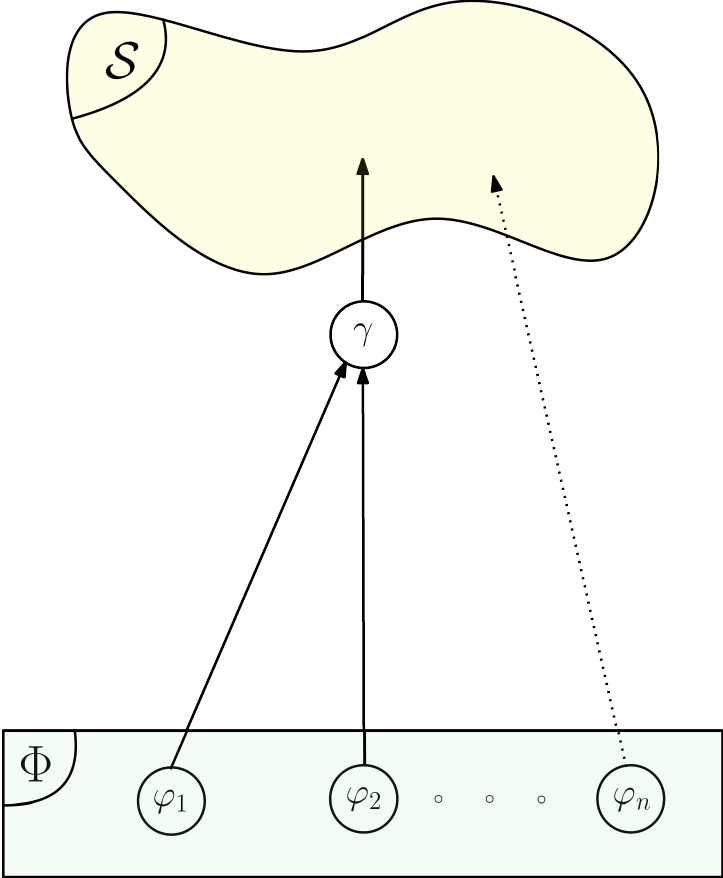}
       \caption{}
        \label{fig:bridge4}
    \end{subfigure}
    \hfill
    \caption{General graphs for the \textit{Bridge} operators. (a) Connection to two STL variables. (b) Connection to an STL variable on the left and a predicate on the right. (c) Connection to a predicate on the left and an STL variable on the right. (d) Direct connection to two predicates. The clouds in blue and yellow represent the rest of the graph. $\mathcal{S}$ is the STL variable space while $\Phi$ represents the root predicate space.}\label{fig:bridges}
\end{figure}

General graphs containing a \textit{Bridge} type operator with STL variable vector $\boldsymbol{\gamma}$ are show in Fig. \ref{fig:bridges}.

The equations governing the transition at the \textit{Bridge} operator level are the following for the four cases represented in Fig. \ref{fig:bridges}:
\begin{equation}
    \begin{matrix}
    \gamma_{k} =\frac{\partial \chi}{\partial \alpha_{k}}\alpha_{k} + \frac{\partial \chi}{\partial \beta_k}\beta_k + R(\bar{\alpha}_k, \bar{\beta}_k)& \text{if} & (\mathcal{P^\gamma_\mathcal{L}} \in \mathcal{S}) \wedge  (\mathcal{P^\gamma_\mathcal{R}} \in \mathcal{S}),\\
    \end{matrix}
\end{equation}
\begin{equation}
    \begin{matrix}
    \gamma_{k} =\frac{\partial \chi}{\partial \alpha_{k}}\alpha_{k} + \frac{\partial \chi}{\partial \rho^{\mathcal{R}}_{k}}\frac{\partial \rho^{\mathcal{R}}_{k}}{\partial \boldsymbol{x}_{\mathcal{R}_{k}}}\boldsymbol{x}_{\mathcal{R}_{k}} + R(\bar{\alpha}_k, \bar{\boldsymbol{x}}_{\mathcal{R}_k})& \text{if} & (\mathcal{P^\gamma_\mathcal{L}} \in \mathcal{S}) \wedge  (\mathcal{P^\gamma_\mathcal{R}} \in \Phi)
    \end{matrix},
\end{equation}
\begin{equation}
    \begin{matrix}
    \gamma_{k} =\frac{\partial \chi}{\partial \rho^{\mathcal{L}}_{k}}\frac{\partial \rho^{\mathcal{L}}_{k}}{\partial \boldsymbol{x}_{\mathcal{L}_{k}}}\boldsymbol{x}_{\mathcal{L}_k} + \frac{\partial \chi}{\partial \beta_k}\beta_k + R(\bar{\boldsymbol{x}}_{\mathcal{L}_k}, \bar{\beta}_k)& \text{if} &  (\mathcal{P^\gamma_\mathcal{L}} \in \Phi) \wedge  (\mathcal{P^\gamma_\mathcal{R}} \in \mathcal{S})
    \end{matrix},
\end{equation}
\begin{equation}
    \begin{matrix}
    \gamma_{k} =\frac{\partial \chi}{\partial \rho^{\mathcal{L}}_{k}}\frac{\partial \rho^{\mathcal{L}}_{k}}{\partial \boldsymbol{x}_{\mathcal{L}_{k}}}\boldsymbol{x}_{\mathcal{L}_{k}} + \frac{\partial \chi}{\partial \rho^{\mathcal{R}}_{k}}\frac{\partial \rho^{\mathcal{R}}_{k}}{\partial \boldsymbol{x}_{\mathcal{R}_{k}}}\boldsymbol{x}_{\mathcal{R}_{k}} + R(\bar{\boldsymbol{x}}_{\mathcal{L}_k}, \bar{\boldsymbol{x}}_{\mathcal{R}_k})& \text{if} & (\mathcal{P^\gamma_\mathcal{L}} \in \Phi) \wedge  (\mathcal{P^\gamma_\mathcal{R}} \in \Phi)
    \end{matrix}.
\end{equation}

In the above expressions  
\begin{itemize}
    \item $\mathcal{P^\gamma_\mathcal{L}}$ (respectively $\mathcal{P^\gamma_\mathcal{R}}$) is the space of the left (respectively right) parent node associated to STL variable $\boldsymbol{\gamma}$; 
    \item $\boldsymbol{x}_{\mathcal{L}}$ (respectively $\boldsymbol{x}_{\mathcal{R}}$) is the state vector of the left (respectively right) predicate;
    \item $\rho^\mathcal{L}$ (respectively $\rho^\mathcal{R}$) is the margin function of the left (respectively right) STL expression.
\end{itemize}

\subsection{General Nesting}

\subsubsection{Backward Propagation}\label{sec:backward-prop}
\begin{figure}[H]
    \centering
    \begin{subfigure}[b]{0.45\textwidth}
        \centering
        \includegraphics[width=1\linewidth]{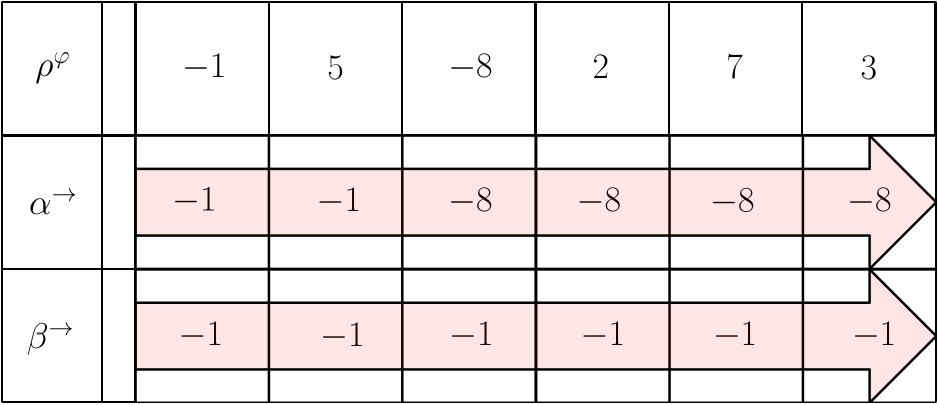}
       \caption{}
        \label{fig:table_forward}
    \end{subfigure}
    \hfill
    \begin{subfigure}[b]{0.45\textwidth}
        \includegraphics[width=1\linewidth]{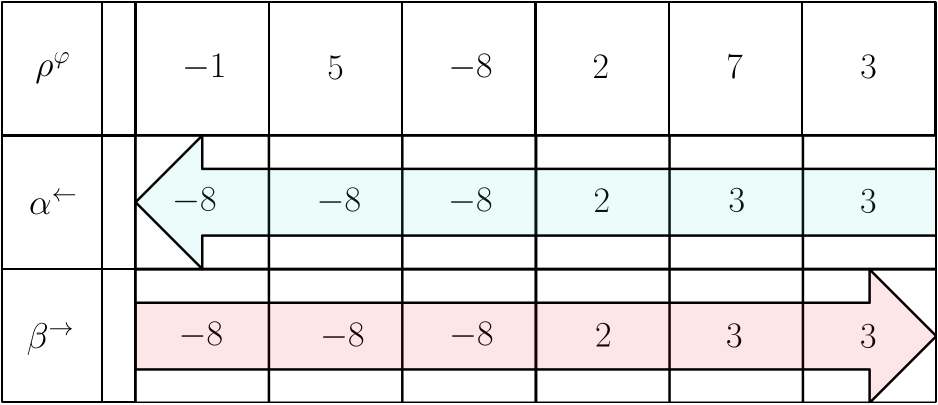}
        \caption{}
        \label{fig:table_backward}
    \end{subfigure}
    \caption{Nesting operator example of Eq. \eqref{eq:ex nesting}. (a) Only forward propagating case. (b) First backward propagating case.}\label{fig:table_forward_backward_1}
\end{figure}

A very important aspect of general nesting is how the margin function is transmitted from each operator to the next one above it. In the previous section, when dealing with a single operator, a forward propagation of the STL variables was performed (i.e., $\alpha_{k+1} = \chi(\alpha_k, \rho^\varphi_{k+1}))$. In the general case, with nesting operators, backward propagation plays a major role in preventing loss of information. To visualize the concept, let us consider the following nesting example:
\begin{equation}\label{eq:ex nesting}
    \underset{\beta}{\diamond} (\underset{\alpha}{\Box} \varphi),   
\end{equation}
represented in Fig. \ref{fig:table_forward_backward_1}, where the notations $\alpha$ and $\beta$ are the STL variables associated to their respective operators. 
To illustrate the concept, let assume that the example result is \textit{True} during the last two time steps, i.e., margin function $\rho^\varphi$ is eventually always positive for these steps. 
As always, the STL variable at the last time step must be positive. One notices that the backward propagation of the margins preserves the Boolean conclusion (i.e., the property is \textit{True}), information being lost at the first negative value in case of forward propagating due to the nesting of $min$ and $max$ functions. 

Another more complex nesting will be presented in Sec. \ref{sec:Until} to build up intuition before presenting the general rules in Sec. \ref{sec:general_formulation}.\\

\subsubsection{\textit{Until} Operator}\label{sec:Until}

\textit{Until} operator, represented with notation $\varphi\mathcal{U}\psi$, and as presented in \cite{scvx_stl}, is an operator that checks if a predicate $\varphi$ is always \textit{True} when another predicate $\psi$ triggered \textit{True}. When taking a closer look at it, $\mathcal{U}$ can be constructed as a nesting of basic operators and written as follows: 
\begin{equation}
    \varphi\mathcal{U}\psi \equiv \diamond (\psi \wedge \Box \varphi).
\end{equation}

On the other hand, when one reads the sentence \textit{When $\psi$ would trigger \textit{True}, $\varphi$ had to have always been \textit{True} before}, a notion of the history of $\varphi$ appears. Therefore, one can state that $\mathcal{U}$ is specifically interested in the history of $\varphi$. For this reason, the \textit{Until} operator could be constructed also as follows
\begin{equation}
    \varphi\mathcal{U}\psi \equiv \diamond (\psi \wedge \Box^\mathcal{B} \varphi),
\end{equation} 
New notation are introduced to take into account temporal direction in the property validity and in the propagation
\begin{itemize}
\item $\Box^\mathcal{B}$ indicates \textit{Always before};
\item $\Box^\mathcal{A}$ indicates \textit{Always after};
\item $\boldsymbol{\alpha}^\rightarrow$ indicates that STL variable $\boldsymbol{\alpha}$ is propagated forward because interested in the past (case of all of what has been done in this paper so far as well as in \cite{scvx_stl});
\item $\boldsymbol{\alpha}^\leftarrow$ indicates that STL variable $\boldsymbol{\alpha}$ is propagated backward because interested in the future.
\end{itemize}

\begin{figure}[h]
        \centering
        \includegraphics[width=0.43\linewidth]{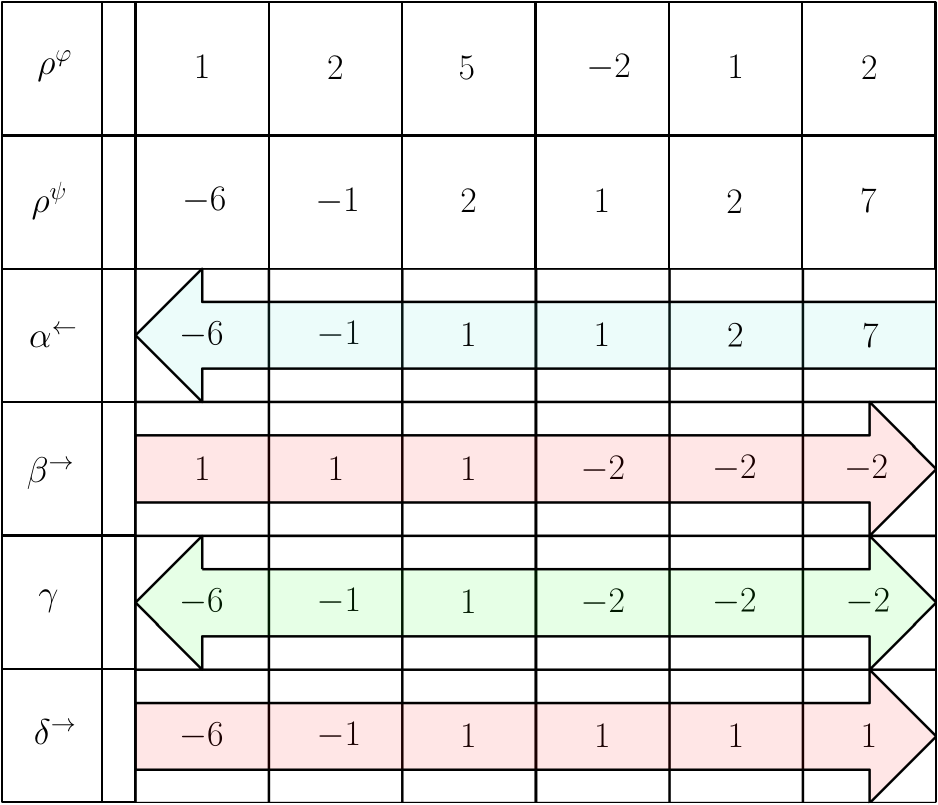}
                \caption{Nesting operator example of Eq. \eqref{nest opr ex}.}
        \label{fig:table2_backward}
    \end{figure}

Another example using \textit{Until} operator is depicted in Fig. \ref{fig:table2_backward} and represented as follows:
\begin{equation}
    \varphi\mathcal{U}(\Box^\mathcal{A}\psi). 
\end{equation} 
This nesting operator is equivalent to the following expression highlighting the propagation direction of the STL variables:
\begin{equation}\label{nest opr ex}
    \underset{\delta^\rightarrow}{\diamond} (\underset{\alpha^\leftarrow}{\Box^\mathcal{A}} \psi) \underset{\gamma}{\wedge} (\underset{\beta^\rightarrow}{\Box^\mathcal{B}} \varphi). 
\end{equation} 
It is recalled that the naming conventions of STL variables follows the order of the Greek alphabet the more external or high in the tree the operator is. Here, the sign of $\delta_N$ respects the given trajectory. 

\subsubsection{General Formulation}\label{sec:general_formulation}

To built a general formulation, about the propagation the following rules can be considered on its direction:
\begin{itemize}
    \item If interested in what happened before when the operator triggers, perform a forward propagation;
    \item If interested in the future when the operator triggers, perform a backward propagation;
    \item Bridges are not affected by propagation direction;
    \item The last node is computed using a forward propagation in any case.
\end{itemize}

\begin{figure}[ht]
    \centering
    \includegraphics[width=0.6\linewidth]{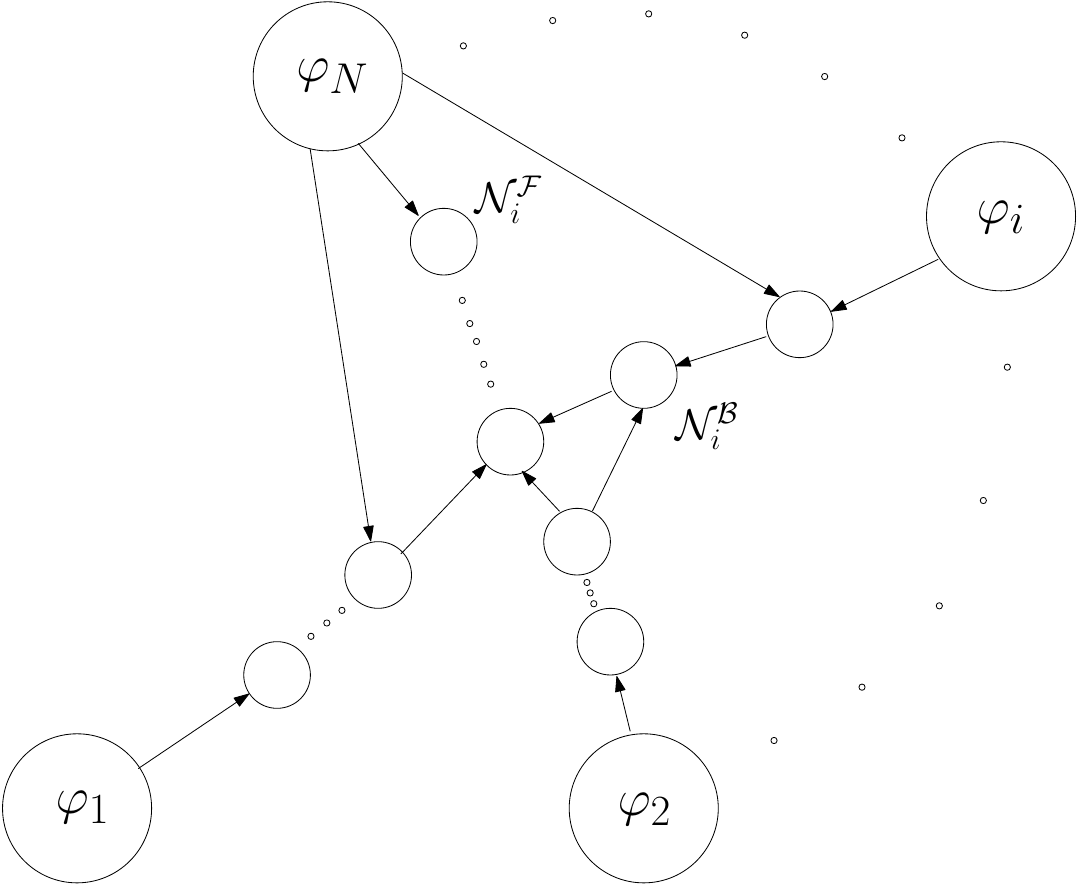}
    \caption{General form of an STL graph.}
    \label{fig:generalization}
\end{figure}

Now that the building blocks are generalized in any position of the tree, full generalization can take place as depicted in Fig. \ref{fig:generalization}. 
\\
Let us define $\mathcal{N}_i$ the generic operator (or node) of index $i$, $\mathcal{B} =\{\wedge, \lor\}$ the \textit{Bridge} operator set, $\mathcal{F} = \{\lozenge, \square\}$ the \textit{Flow} operator set, $\varphi_i$ the root predicates. One then gets
\\
\begin{equation}
\begin{cases}
  \alpha_{k}^{\mathcal{N}_i^\mathcal{B}} = A_{k}^{\mathcal{P}_\mathcal{L}\left(\mathcal{N}_i^\mathcal{B}\right)}\alpha_{k}^{\mathcal{P}_\mathcal{L}\left(\mathcal{N}_i^\mathcal{B}\right)} + A_{k}^{\mathcal{P}_\mathcal{R}\left(\mathcal{N}_i^\mathcal{B}\right)}\alpha_{k}^{\mathcal{P}_\mathcal{R}\left(\mathcal{N}_i^\mathcal{B}\right)}\\
    \alpha_{k-1}^{\mathcal{N}_i^\mathcal{F}} = A_k^{\mathcal{N}_i^\mathcal{F}}\alpha_k^{\mathcal{N}_i^\mathcal{F}} + A_{k-1}^{\mathcal{P}\left(\mathcal{N}_i^\mathcal{F}\right)} \alpha_{k-1}^{\mathcal{P}\left(\mathcal{N}_i^\mathcal{F}\right)} & \mathrm{if} \;\;\; \left(\mathcal{C}\left(\mathcal{N}_i\right) \neq \varnothing\right) \lor \left(\mathcal{T}(\mathcal{N}_i^\mathcal{F}) = \mathcal{A}\right)\\
    \alpha_{k+1}^{\mathcal{N}_i^\mathcal{F}} = A_k^{\mathcal{N}_i^\mathcal{F}}\alpha_k^{\mathcal{N}_i^\mathcal{F}} + A_{k+1}^{\mathcal{P}\left(\mathcal{N}_i^\mathcal{F}\right)} \alpha_{k+1}^{\mathcal{P}\left(\mathcal{N}_i^\mathcal{F}\right)} & \mathrm{if} \;\;\; \left(\mathcal{C}\left(\mathcal{N}_i\right) = \varnothing \right) \lor \left(\mathcal{T}(\mathcal{N}_i^\mathcal{F}) = \mathcal{B}\right)\\
\end{cases},
\end{equation}
\\
where 
\begin{itemize}
    \item  $\alpha_{k}^{\mathcal{N}_i^\mathcal{F}}$ and $\alpha_{k}^{\mathcal{N}_i^\mathcal{B}}$ are respectively the STL variables associated to the \textit{Flow} and \textit{Bridge} operators at position $i$ and time step $k$;
    \item $A^{\mathcal{N}_i}_k$ are the Jacobian matrices associated to each operator $\mathcal{N}_i$ in the tree at time step $k$;
    \item $\mathcal{P}(\mathcal{N}_i)$ (respectively $\mathcal{C}(\mathcal{N}_i)$) are the parent operator (respectively child operator) of node $\mathcal{N}_i$ (it must be noted that for a predicate, by definition at the bottom of the tree, this coincides exactly with the margin function of the states throughout the trajectory);
    \item $\mathcal{T}(\mathcal{N}_i^\mathcal{F})$ the temporal interest of \textit{Flow} type node $\mathcal{N}_i^\mathcal{F}$, being equal to $\mathcal{A}$ for \textit{After} and $\mathcal{B}$ for \textit{Before};
    \item $\varnothing$ the empty set.
\end{itemize}

Now, since all the equations are linearized, and since an STL variable solely depends on its parents, it is possible to write the entire constraint using block matrices and vectors. As example, the following nested operator is considered:
\begin{equation}\label{new nest opr ex}
(\varphi \underset{\alpha}{\wedge} \psi) \underset{\gamma^\rightarrow, \delta, \epsilon^\rightarrow}{\mathcal{U}} (\underset{\beta^\leftarrow}{\Box} \zeta).  
\end{equation}

    \begin{figure}[ht]
        \centering
        \includegraphics[width=0.30\linewidth]{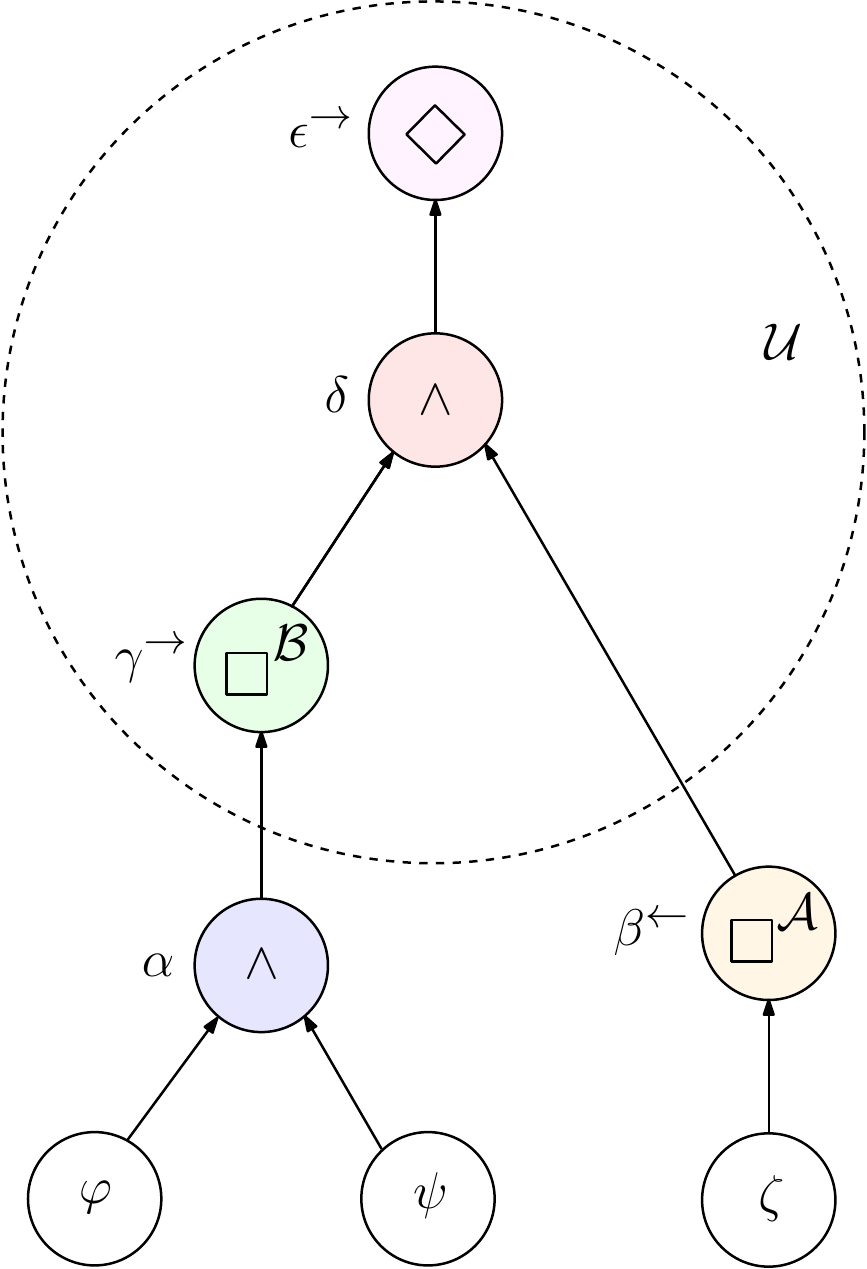}
        \caption{Graph representation of the nested operator given by Eq. \eqref{new nest opr ex}.}
        \label{fig:example_graph}
    \end{figure}
In reference to Fig. \ref{fig:example_graph}, the entire state-transition matrix has the following structure

\begin{equation}
    \begin{bmatrix}
        \boldsymbol{A^{\alpha}_{x}} & \boldsymbol{A^{\alpha}_{\alpha}} & \boldsymbol{0_{N\times N}}& \boldsymbol{0_{N\times N}} & \boldsymbol{0_{N\times N}}& \boldsymbol{0_{N\times N}}\\
        \boldsymbol{A^{\beta}_{x}} &
        \boldsymbol{0_{N \times N}} & 
        \boldsymbol{A^{\beta}_{\beta}} & 
        \boldsymbol{0_{N\times N}}& \boldsymbol{0_{N\times N}}& \boldsymbol{0_{N\times N}}\\
        \boldsymbol{0_{N\times N}} & \boldsymbol{A^{\gamma}_{\alpha}}& \boldsymbol{0_{N\times N}} &
        \boldsymbol{A^{\gamma}_{\gamma}} & \boldsymbol{0_{N\times N}}& \boldsymbol{0_{N\times N}}\\
        \boldsymbol{0_{N\times N}} & \boldsymbol{0_{N\times N}}  & \boldsymbol{A^{\delta}_{\beta}} &
        \boldsymbol{A^{\delta}_{\gamma}}  & \boldsymbol{A^{\delta}_{\delta}}& \boldsymbol{0_{N\times N}}\\
        \boldsymbol{0_{N\times N}} &  \boldsymbol{0_{N\times N}}& \boldsymbol{0_{N\times N}} &
        \boldsymbol{0_{N\times N}}& \boldsymbol{A^{\epsilon}_{\delta}} & \boldsymbol{A^{\epsilon}_{\epsilon}} 
    \end{bmatrix}
    \begin{bmatrix}
        \boldsymbol{x}\\ \boldsymbol{\alpha}\\ \boldsymbol{\beta} \\ \boldsymbol{\gamma} \\ \boldsymbol{\delta} \\ \boldsymbol{\epsilon}
    \end{bmatrix}
     =
     \begin{bmatrix}
         \boldsymbol{b^\alpha} \\
         \boldsymbol{b^\beta}\\
         \boldsymbol{b^\gamma}\\
         \boldsymbol{b^\delta}\\
         \boldsymbol{b^\epsilon}

     \end{bmatrix}.
\end{equation} with $\boldsymbol{A}_\bullet^\bullet$ and $\boldsymbol{b}_\bullet^\bullet$ defined following the theory presented in Sections \ref{sec:isolated_operators} and \ref{sec:nesting}. Finally, on top of this block-matrix equality, the inequality $\epsilon_N \geq 0$ is added.\\

\begin{figure}[H]
    \centering
    \begin{subfigure}[b]{0.45\textwidth}
        \centering
        \includegraphics[width=1\linewidth]{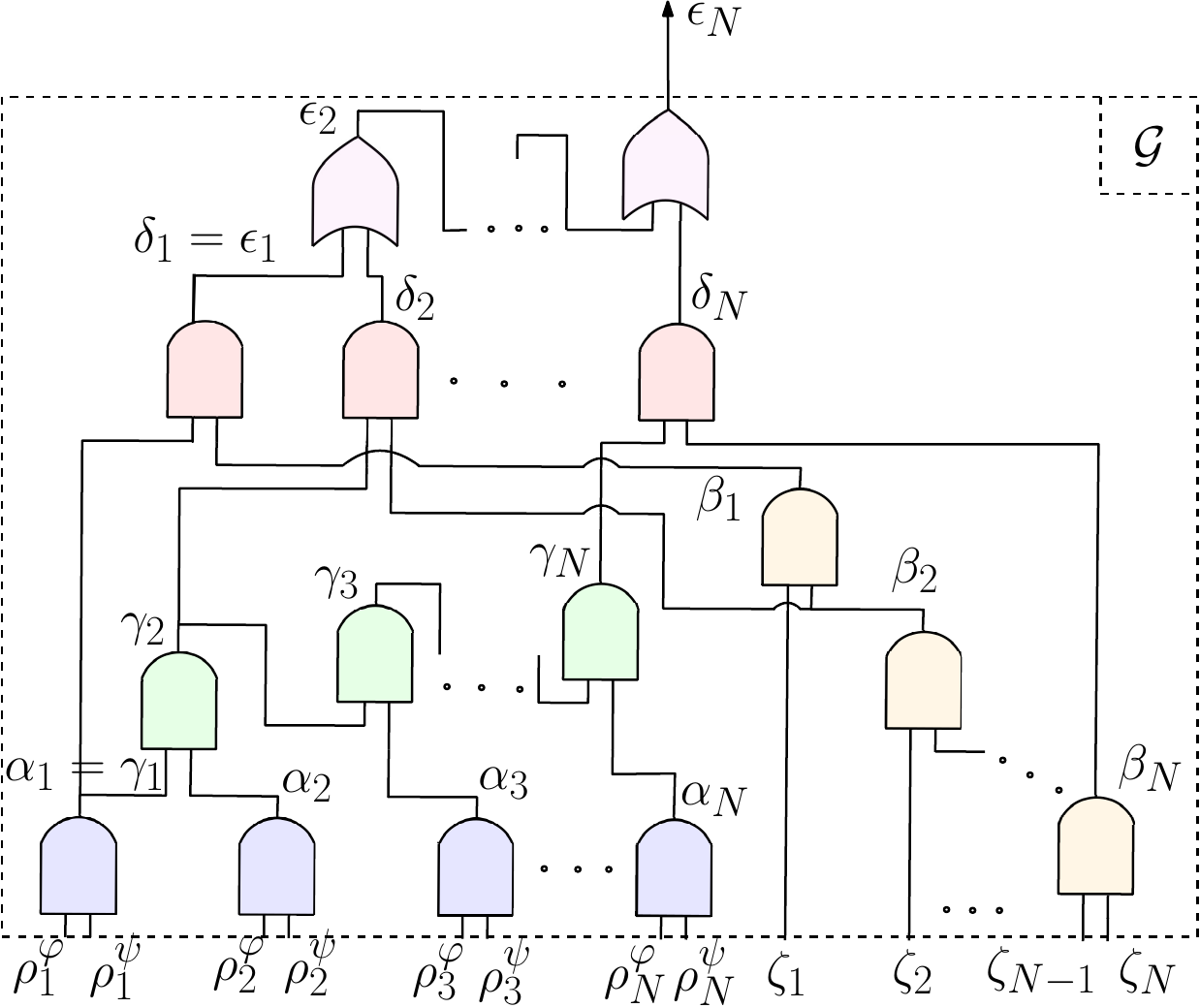}
       \caption{}
        \label{fig:circuit}
    \end{subfigure}
    \hfill
    \begin{subfigure}[b]{0.45\textwidth}
        \includegraphics[width=1\linewidth]{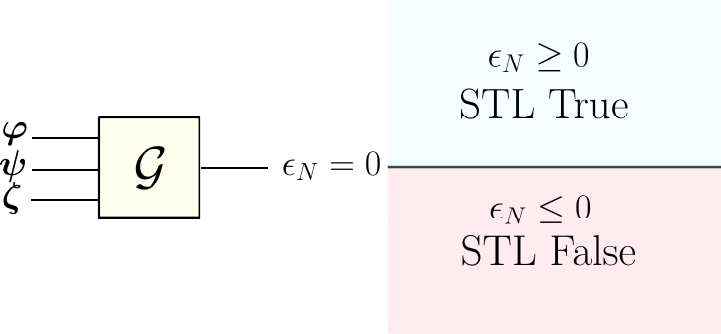}
        \caption{}
        \label{fig:circuit_condensed}
    \end{subfigure}
    \caption{Circuit graph representation of the nested operator given in Eq. \eqref{new nest opr ex}. (a) Electronic analogy. (b) Input-output equivalence, condensed form.}\label{fig:circuits}
\end{figure}

Graph representation of Fig. \ref{fig:example_graph} can also be seen from the point of view of an electronic circuit, as shown in Fig. \ref{fig:circuit} and Fig. \ref{fig:circuit_condensed} (condensed form). Both representations are equivalent. From several predicates, the output is a number. Its sign defines the validity of the whole constraint. 
    

\section{Derived Operators}
This section presents three more operators which can be expressed in function of the building blocks already defined. They can serve as a more visual way of modeling a problem by regrouping operators together.

\subsection{\textit{Implies} ($\implies$)}
\textit{Implies} ($\implies$) can be considered as a bridge operator linking two predicates $\varphi$ and $\psi$. 
Using the STL grammar the \textit{Implies} operator can be expressed as follows:
\begin{equation}\label{eq implies}
    \varphi \implies \psi \equiv (\neg \varphi) \lor \psi.
\end{equation}
One should just be careful in what they mean when using \textit{Implies}. 
Eq. (\ref{eq implies}) logically means either $\varphi$ is \textit{False}, or, $\varphi$ and $\psi$ are \textit{True} at the same time. This detail can be of great use to simplify graphs.

\subsection{\textit{If and Only If} ($\iff$)}
\textit{If and Only If} ($\iff$) can be considered as a double implication. Using the STL grammar it can be expressed as follows:
\begin{equation}
    \varphi \iff \psi \equiv (\varphi \implies \phi) \wedge (\psi \implies \varphi) \equiv (\neg \varphi \lor \psi) \wedge (\neg \psi \lor \varphi).
\end{equation}

\subsection{\textit{Exclusive Or} ($\oplus$)}
Operator \textit{Exclusive Or} ($\oplus$) can be considered as a variant of the inclusive disjunction. It is also the negation of operator \textit{If and Only If} for the De-Morgan's law. Therefore, operator $\oplus$ can be defined either as a nesting of disjunctions and conjunctions, or as a negation of the equivalence:
\begin{equation}
    \varphi \oplus \psi \equiv [\varphi \wedge (\neg \psi)] \lor [(\neg \varphi \wedge \psi)] \equiv \varphi \neg(\iff) \psi.
\end{equation}

\subsection{Other Extensions}
Defining new operators can be made using previously defined ones in the same philosophy as for \textit{Until}. Other ideas could be neither nor, exclusive negation or, negative and, or maybe more complex operations like $\varphi \mathcal{U}(\psi \wedge \Box (\diamond \zeta))$ or $\diamond (\Box (\diamond \varphi))$. Extensions can be created and then used as a block by itself to form more complex logic, hence manipulating super-operators and getting closer to real-world operations.

\section{Applications}\label{sec:applications}
This section presents two examples of STL constraints used in the context of the satellite trajectory optimization. The modified successive convexification scheme is presented in detail in Section \ref{sec:modif-scvx}, implemented for two satellite rendezvous mission applications and discussed in Sections \ref{sec:first_app} and \ref{sec:second_app}.

\subsection{Problem Setup}
As general definition of a satellite rendezvous mission, one satellite (named as \textit{chaser}) has to dock to another one (named as \textit{target}). Here, the \textit{target} is in a circular geostationary orbit around the Earth. The \textit{chaser} has a mass of 1000 kg and it is initially at coordinates [-2000, 500, 200] m in an orbital frame centered in the \textit{target} position and with a null relative velocity. The problem is of minimum-fuel optimization over a finite horizon of 5000 s. The maximum allowed thrust of the \textit{chaser} is 1 N. The simulations were performed using Matlab R2021a and the Mosek solver \cite{mosek}. The initial trajectory was obtained with zero thrust (i.e., it is a free drift propagation). All the details of this problem statement can be found in \cite{claudet2023}.

\subsection{Modified Successive Convexification Algorithm}\label{sec:modif-scvx}
\begin{figure}[ht]
        \centering
        \includegraphics[width=0.6\linewidth]{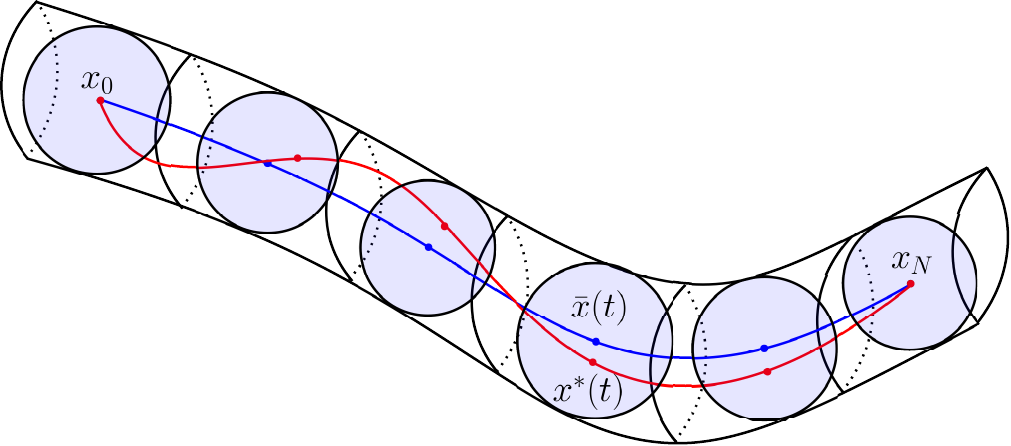}
        \caption{Trust region. The reference trajectory $\bar{\boldsymbol{x}}(t)$ is the blue line while the new optimal solution $\boldsymbol{x}^*(t)$) is the red one. The light blue circles around the blue dots represent the trust regions around which each point of the next solution must lie in.}
        \label{fig:trust_region}
    \end{figure}

Since the problem is in a convex form, the choice was made to use the powerful Successive Convexification (SCvx) framework which extends Sequential Convex Programming (SQP) using trust regions (shown on Fig. \ref{fig:trust_region}) for improved performance \cite{scvx}. The algorithm, which is a modified version of the one in \cite{claudet2023} is presented in Algorithm \ref{alg:SCvx}.\\

To integrate STL constraint in the trajectory optimization problem, the following hypotheses have been assumed.

\begin{figure}[h]
    \centering
    \begin{subfigure}[b]{0.40\textwidth}
        \centering
        \includegraphics[width=1\linewidth]{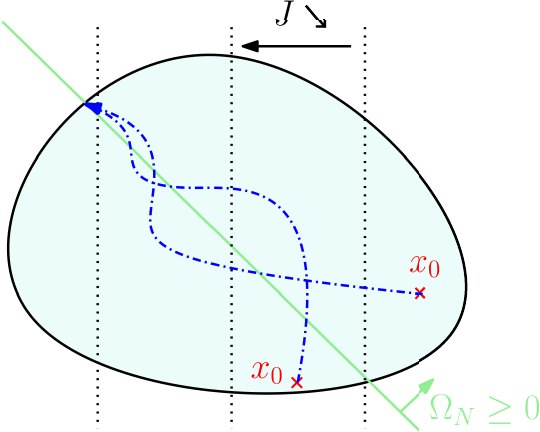}
       \caption{}
        \label{fig:STL_convergence_1}
    \end{subfigure}
    \hfill
    \begin{subfigure}[b]{0.40\textwidth}
        \includegraphics[width=1\linewidth]{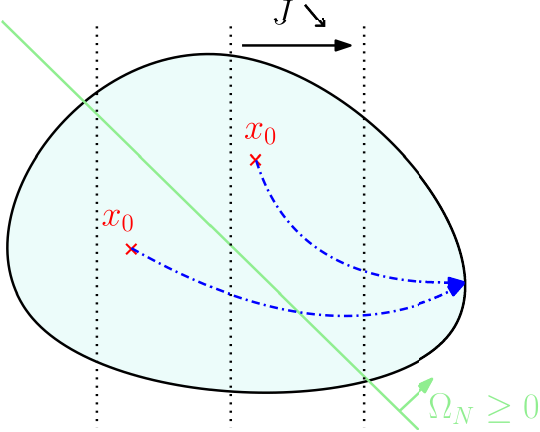}
        \caption{}
        \label{fig:STL_convergence_2}
    \end{subfigure}
    \caption{Theoretical behavior of the solution when adding the STL constraint. Value $x_0$ represents the initial guess, the cost lines are in dotted black, the successive iterations are shown in blue, and the STL plane in green.}\label{fig:STL_congergence_entire}
\end{figure}

\begin{itemize}
    \item As explained in Section \ref{subsubsec:linearization}, the constraints are functions of the reference trajectory appearing in the zero order term and in the evaluations of the gradients of the linearization. Therefore, a new optimal solution is only optimal with respect to the reference trajectory. To make sure that the new optimal solution really satisfies the STL constraint, the idea is to compute the margin of the STL constraint with the new updated optimal values. This is not a constraint anymore since optimization is already done for this iteration, but a model-checking step. Note that, as the solution converges, this check is less and less important (indeed the reference becomes closer to the new optimal trajectory). 
    \item The cost function to minimize is modified from
    \begin{equation}
    J = \mathcal{E}
\end{equation}
to 
   \begin{equation}\label{cost}
    J_{STL} = \mathcal{E} + w_{\Omega}\Omega_N.
\end{equation}
With respect to the energy $\mathcal{E}$, the additional second term is function of the STL variable $\Omega_N$ associated to the last time step of the last child at the top of the graph, and of its associated weight $w_\Omega$. As a design choice for the cost function, it was decided for $\Omega_N$ to increase its weight $w_\Omega$ negatively until a first solution satisfying the STL constraint is found. Indeed, to minimize $J_{STL}$ with $w_\Omega < 0$, $\Omega_N$ has to be positive (i.e., constraint satisfied). After that, as long as the STL constraint stays \textit{True} ($\Omega_N > 0)$, $w_\Omega$ is increased, getting closer to zero from negative value, in orde to give more importance to the other objectives as the energy minimization. When $\Omega_N$ becomes negative, the STL constraint becomes \textit{False} by construction and more importance is again given to the STL constraint, in order to converge to the STL plane $\{\Omega_N = 0\}$ while minimizing the energy to find the optimal solution (see Fig. \ref{fig:STL_congergence_entire} for a graphical representation).
    \item Moreover, it is also possible to play on the trust radius (see Fig. \ref{fig:trust_region}). The radius is increased until a solution satisfying the STL constraint is found, then is shrinks down to help convergence.
\end{itemize}
The modifications to the Successive Convexification algorithm present in \cite{claudet2023} are in red within Algorithm \ref{alg:SCvx}'s description:
  \begin{algorithm}
            \caption{Modified Successive Convexification (SCvx)}\label{alg:SCvx}
                \begin{algorithmic}[5]
                 \small{
                    \State Initialization
                    \State Scaling of variables
                    \While{Not Accurate or Not Converged}
                            \State Linearization 
                            \State Discretization
                            \State Constraints update
                            \State Solving of the convexified sub-problem
                            \State Non-linear repropagation
                            \State \textcolor{red}{STL model-checking} 
                        \If{Not Feasible}
                            \State Dilate trust region 
                            \State\textcolor{red}{Update cost weights} 
                        \Else 
                        \If{Not Accurate}
                                \State Update trust region
                                \State\textcolor{red}{Update cost weights} 
                            \Else \If{Not Converged}
                                \State Update trust region
                                \State Update reference
                                 \State\textcolor{red}{Update cost weights} 
                            \EndIf
                        \EndIf
                        \EndIf
                    \EndWhile}
                \end{algorithmic}
            \end{algorithm}                    

The theoretical behavior of the solution when adding the STL is shown in Fig. \ref{fig:STL_congergence_entire}. In the case of a contradiction between satisfying the STL constraint and minimizing the cost (most relevant case, depicted in Fig. \ref{fig:STL_convergence_1}), the iterations will converge to the STL plane $\{\Omega_N = 0\}$, while minimizing the cost. In the other case (depicted in Fig. \ref{fig:STL_convergence_2}), when the natural optimal solution is also fulfilling the STL requirement, the solution converges to the minimal cost, satisfying the STL constraint but not lying on its plane.

\subsection{First Example: Safety Spacing in a Satellite Rendezvous Mission}\label{sec:first_app}
\begin{figure}[H]
        \centering
        \includegraphics[width=1\linewidth]{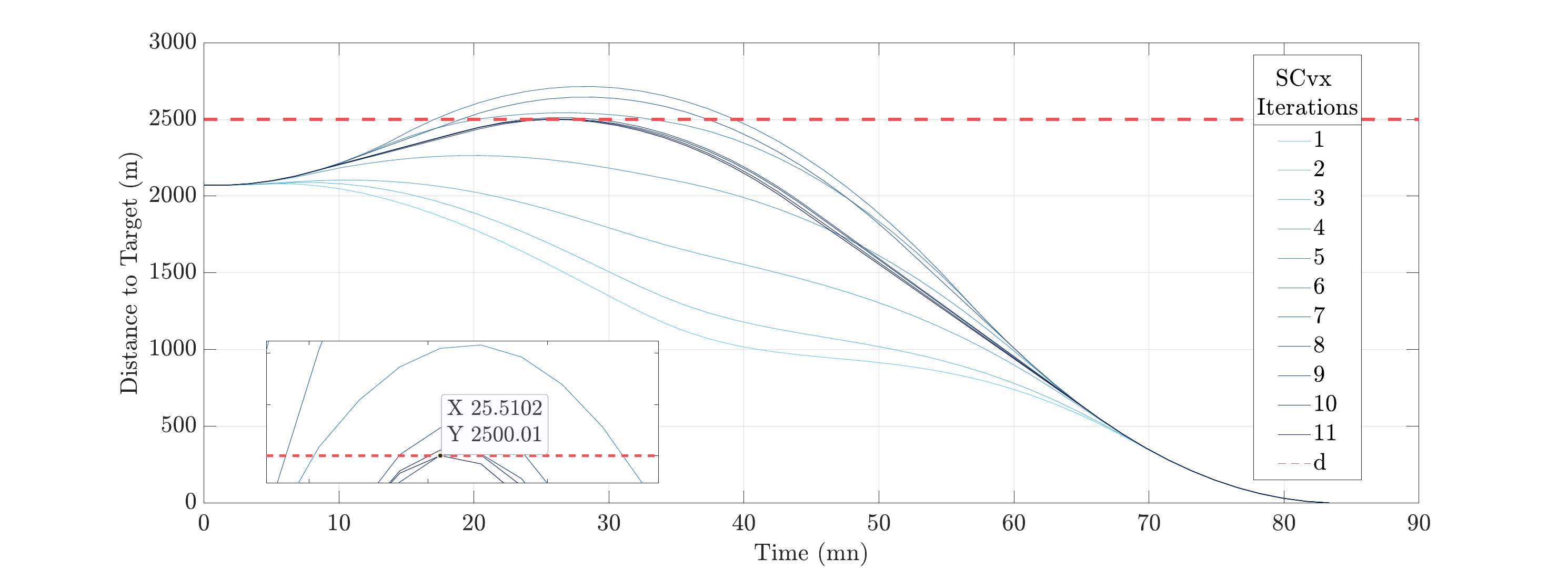}
        \caption{Iterations of the SCvx algorithm to solve the minimum-fuel optimization while imposing a safety spacing of 2500 m (red dashed line) around the \textit{target}.}
        \label{fig:distance_greater_states}
    \end{figure}

\begin{figure}[ht]
     \centering
     \begin{subfigure}[b]{0.49\textwidth}
         \centering
        \includegraphics[width=1\linewidth]{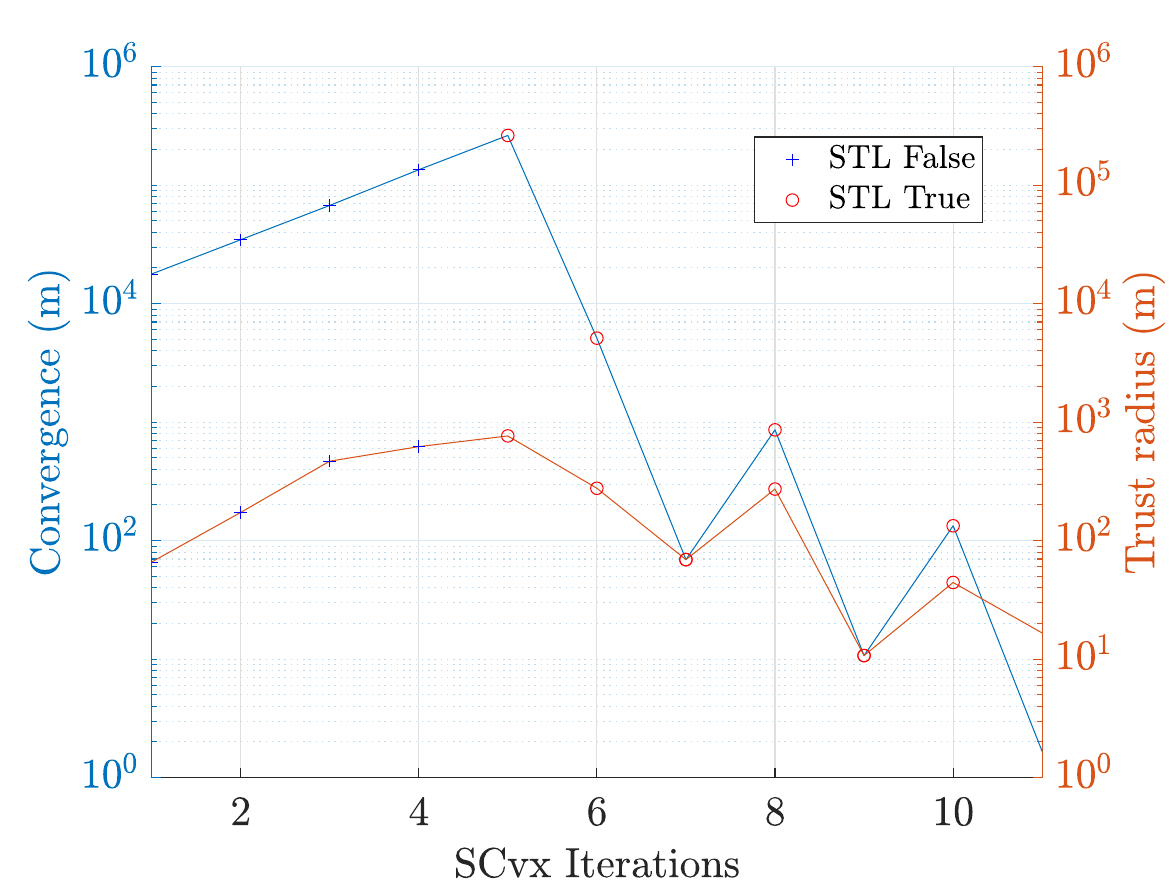}
            \caption{}
        \label{fig:distance_greater_convergence_trust_radius}
     \end{subfigure}
     \hfill
     \begin{subfigure}[b]{0.49\textwidth}
         \centering
        \includegraphics[width=1.\linewidth]{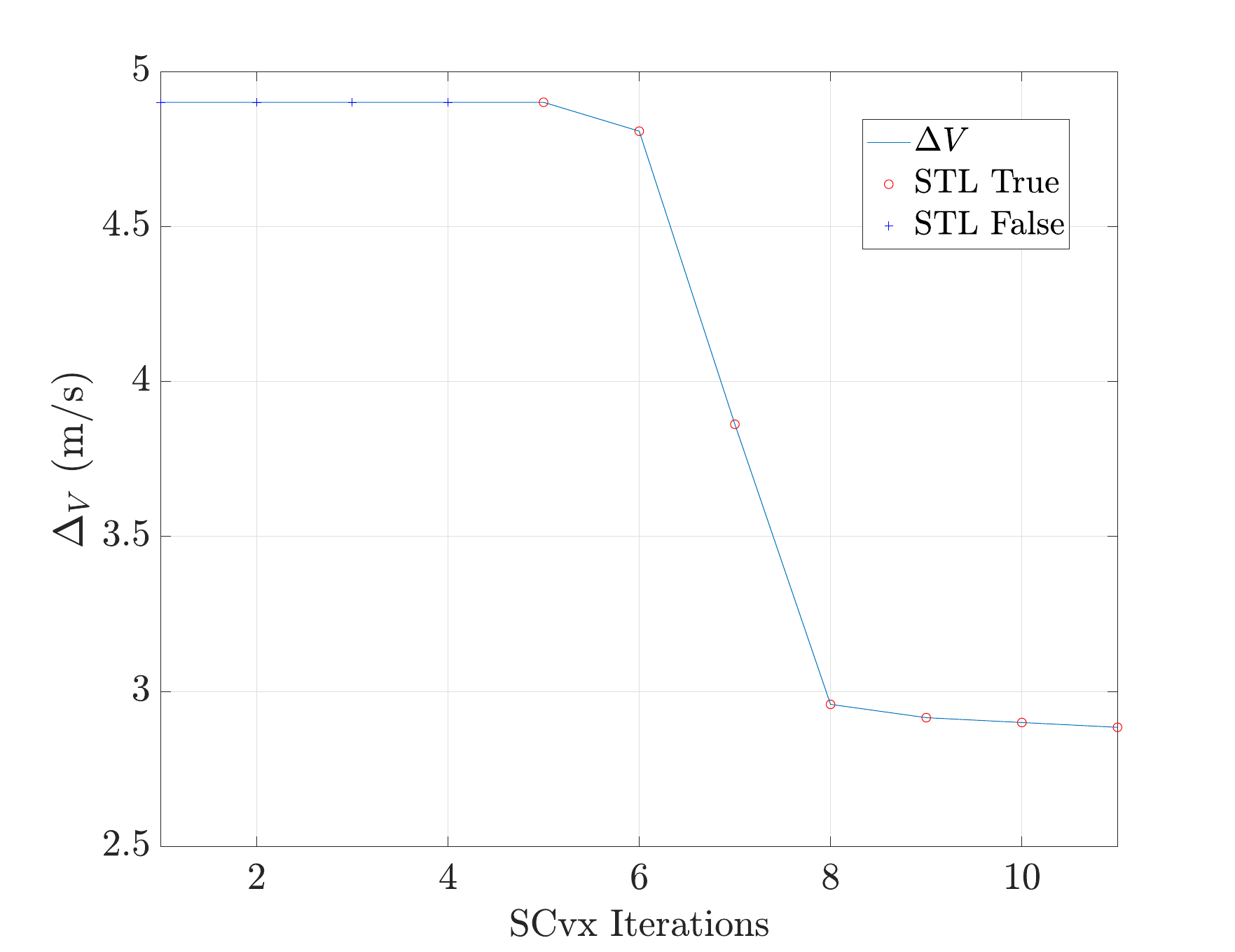}
            \caption{}
        \label{fig:distance_greater_delta_v}
     \end{subfigure}
        \caption{(a) Convergence and trust radius. (b) Energy minimization.}
        \label{fig:convergence_trust_radius_energy}
\end{figure}

In the first example, let us imagine that certain sensors (e.g., GPS or lidar) during the approach experience a failure. In that case, a possible safety maneuver to prevent collision to the \textit{target} could be for the \textit{chaser} to first go outside of a safety sphere around the target before trying to approach and dock it again. While maintaining the same other constraints and objectives, this constraint was modeled as \textit{Eventually get outside of a circle of radius $d$}, or, using the STL grammar, as:
\begin{equation}
    \diamond \lVert \boldsymbol{r}_{ct} \rVert \geq d,
\end{equation} 
where $\boldsymbol{r}_{ct}$ is the relative position vector between the \textit{chaser} and the \textit{target} in the relative orbital frame, and $d$ is the radius of the safety sphere. The iterations of the successive trajectories can be seen in Fig. \ref{fig:distance_greater_states}. As one could guess, the optimal trajectory tangents to the safety sphere before going back and dock. In this context, convergence value $\epsilon_c$ is defined as:

 \begin{equation}
        \epsilon_c = \max_{i=1 \dots N} \left\lVert 
                    \boldsymbol{r^{*(i)}}_{ct}- \boldsymbol{\bar{r}^{(i)}}_{ct}\right \rVert,
    \end{equation}    
while trust radius $TR$ is a standard constraint of the SCvx algorithm, written for each iteration as:

\begin{equation}
                \forall i \in [1, \hdots, N], \;\; \left\lVert 
                    \boldsymbol{r^{*(i)}}_{ct} - \boldsymbol{\bar{r}^{(i)}}_{ct} \right \rVert \leq TR^{(i)},
            \end{equation}
    with $\boldsymbol{r^{*}}_{ct}$ the optimal value of the position and $\boldsymbol{\bar{r}}_{ct}$ its reference value (i.e., the previous optimal solution). Fig. \ref{fig:distance_greater_convergence_trust_radius} shows convergence (with a threshold of $\epsilon_c = 5$ m) and trust radius (no specific threshold) while the $\Delta$V representing the fuel consumption is observed in Fig. \ref{fig:distance_greater_delta_v}. The minimum-fuel optimal solution is found in two main steps. First, a feasible solution is obtained where the STL constraint is satisfied. From that point on, the solution is adjusted to keep satisfying the STL constraint, but this time giving more and more importance to the other objective (i.e., energy minimization). \\


\subsection{Second Example: Slow Approach in a Satellite Rendezvous Mission}\label{sec:second_app}

The second example focuses on a nesting operator using a \textit{Bridge} and a \textit{Flow} operator. The objectives id that between 60\% and 80\% of the transfer time, the relative distance between \textit{target} and \textit{chaser} is always less than $R=500$ m, and that the norm of the thrust vector $\boldsymbol{F}$ is below the half of its maximum value of 1 N. Using the STL grammar, these constraints can be formalized as follows:  
\begin{equation}
    \Box _{[60\%, 80\%]} \{(\lVert \boldsymbol{r}_{ct} \rVert \leq R) \wedge (\lVert \boldsymbol{F} \rVert \leq F_{max}/2) \}.
\end{equation}

\begin{figure}[h!]
     \centering
     \begin{subfigure}[b]{0.45\textwidth}
         \centering
        \includegraphics[width=1\linewidth]{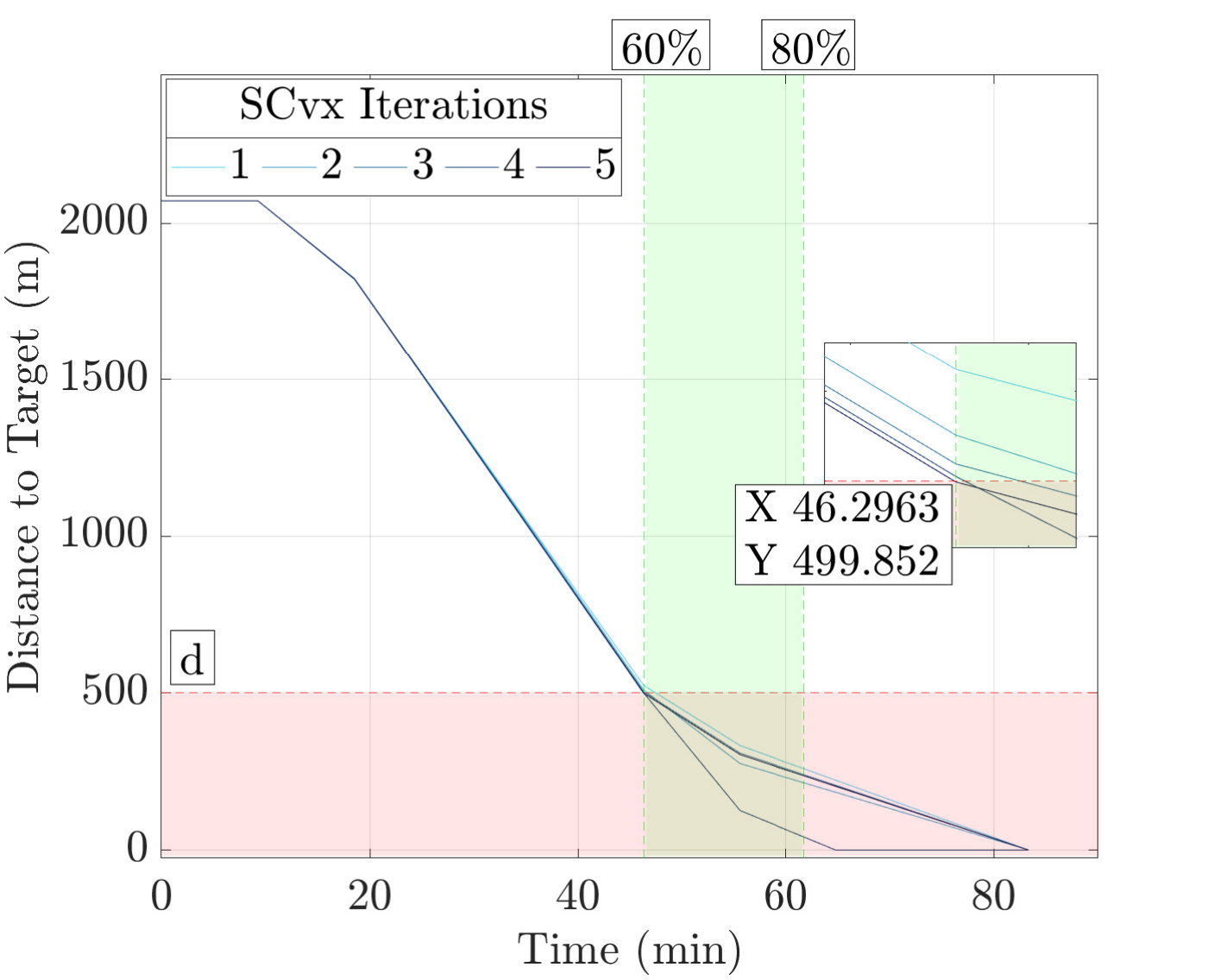}
            \caption{}
        \label{fig:always_thrustmaxover2_and_dlessthanR_at60percent_new}
     \end{subfigure}
     \hfill
     \begin{subfigure}[b]{0.45\textwidth}
         \centering
        \includegraphics[width=1.\linewidth]{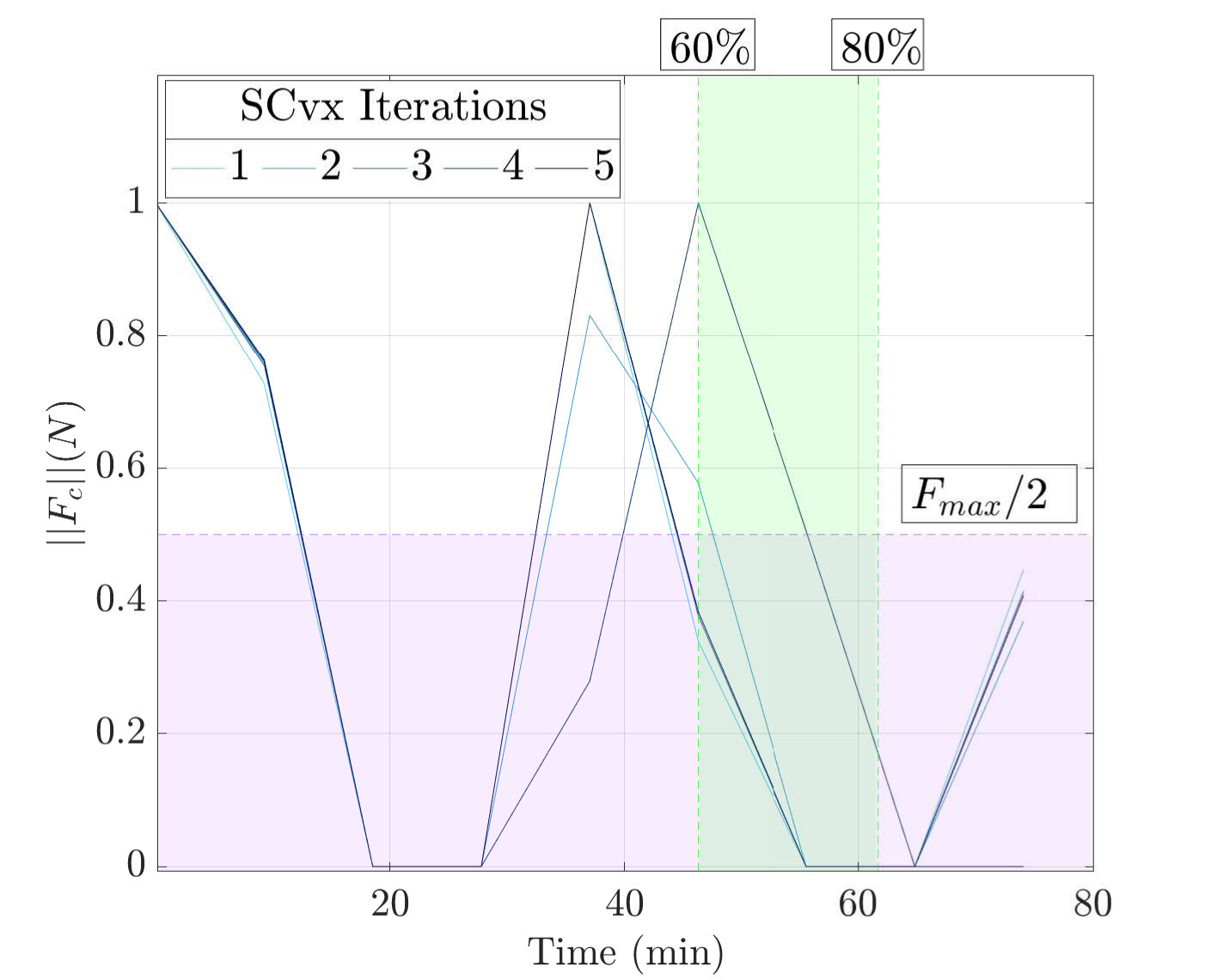}
            \caption{}
        \label{fig:thrust_for_always_thrustmaxover2_and_dlessthanR_at60percent_new}
     \end{subfigure}
        \caption{(a) Distance of the chaser to the target over time. (b) Thrust norm over time.}
        \label{fig:always_thrustmaxover2_and_dlessthanR_at60percent}
\end{figure}

Fig. \ref{fig:always_thrustmaxover2_and_dlessthanR_at60percent_new} shows the distance of the \textit{chaser} to the \textit{target} over time. As one would think, the iterations converge to the threshold at one of the time bounds. The thrust norm over time is represented in Fig. \ref{fig:thrust_for_always_thrustmaxover2_and_dlessthanR_at60percent_new}. The iterations converge meeting the constraint.\\

\begin{figure}[h]
        \centering
        \includegraphics[width=1\linewidth]{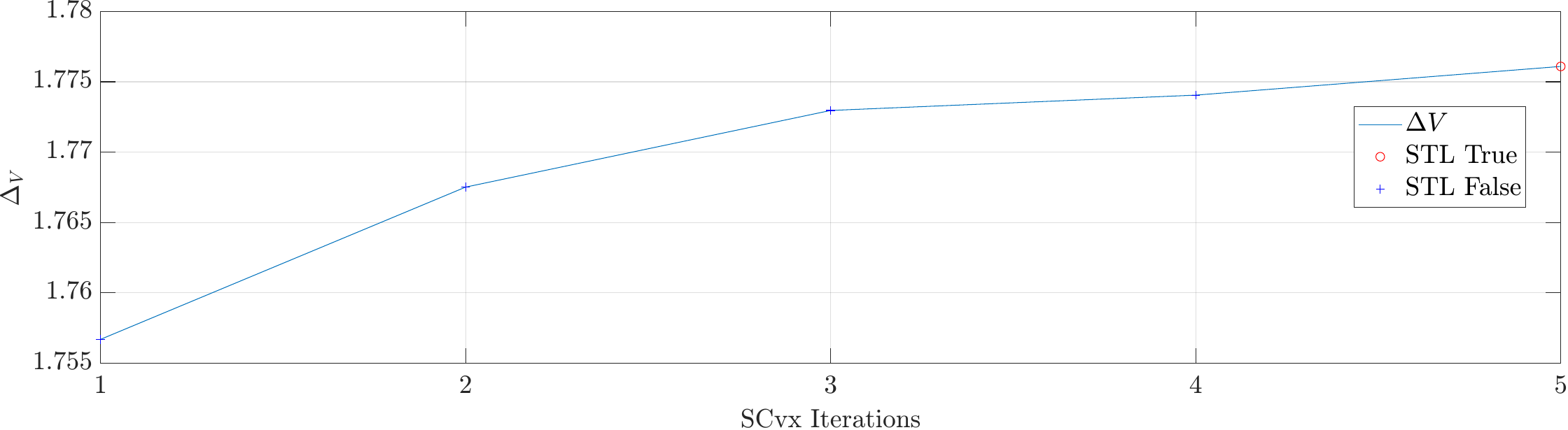}
        \caption{Total fuel consumption at each iteration.}
        \label{fig:DeltaV_for_always_thrustmaxover2_and_dlessthanR_at60percent}
    \end{figure}

Then, in Fig. \ref{fig:DeltaV_for_always_thrustmaxover2_and_dlessthanR_at60percent}, the evolution of the fuel consumption is shown at each iteration. It increases until the STL constraint becomes positive. 

\newpage
\section{Conclusion}\label{sec:conclusion}

In this work, an end-to-end methodology was presented to be able to take into account arbitrarily-complex Signal Temporal Logic constraints in optimization problem solving. It was shown that basic logic operators such as \textit{Or, And, Eventually, Always} can be nested together, discretized and linearized to fit powerful convex frameworks. Moreover, graph-based STL formalism is well suited for the construction of super-operators, which, in turn, ease the understanding and the design of complex high-level mission scenarios and pave the way towards realistic agile operations. Using a modified Successive Convexification scheme, two examples were discussed in the context of the safety in autonomous spacecraft rendezvous missions. The results demonstrated high precision and fast convergence properties. 

\section{Acknowledgment}

This paper is part of a Ph.D. work in progress at ISAE-SUPAERO, funded by the IRT (Technological Research Institute) of Saint-Exupéry, in the context of the RAPTOR project (Robotic and Artificial intelligence Processing Test on Representative target) for satellite rendezvous missions, in collaboration with Thales Alenia Space.

\bibliographystyle{alpha}
\bibliography{sample}

\newcommand{\etalchar}[1]{$^{#1}$}
\begin{thebibliography}{YYW{\etalchar{+}}19}

\bibitem[ApS23]{mosek}
MOSEK ApS.
\newblock {\em The MOSEK optimization toolbox for MATLAB manual. Version
  10.0.}, 2023.

\bibitem[BV04]{boyd2004}
Stephen Boyd and Lieven Vandenberghe.
\newblock {\em Convex Optimization}.
\newblock {Cambridge University Press}, March 2004.

\bibitem[CLL{\etalchar{+}}23]{claudet2023}
Thomas Claudet, Jérémie Labroquère, Damiana Losa, Francesco Sanfedino, and
  Daniel Alazard.
\newblock Successive convexification for on-board scheduling of spacecraft
  rendezvous missions.
\newblock ESA GNC, 2023.

\bibitem[DM10]{donze2010}
Alexandre Donz{\'e} and Oded Maler.
\newblock Robust satisfaction of temporal logic over real-valued signals.
\newblock In Krishnendu Chatterjee and Thomas~A. Henzinger, editors, {\em
  Formal Modeling and Analysis of Timed Systems}, pages 92--106, Berlin,
  Heidelberg, 2010. Springer Berlin Heidelberg.

\bibitem[LAP23]{leung2023backpropagation}
Karen Leung, Nikos Ar{\'e}chiga, and Marco Pavone.
\newblock Backpropagation through signal temporal logic specifications:
  Infusing logical structure into gradient-based methods.
\newblock {\em The International Journal of Robotics Research}, 42(6):356--370,
  2023.

\bibitem[MAGC22]{scvx_stl}
Yuanqi Mao, Behcet Acikmese, Pierre-Loic Garoche, and Alexandre Chapoutot.
\newblock Successive convexification for optimal control with signal temporal
  logic specifications.
\newblock In {\em Proceedings of the 25th ACM International Conference on
  Hybrid Systems: Computation and Control}, HSCC '22, New York, NY, USA, 2022.
  Association for Computing Machinery.

\bibitem[MGFS20]{Ma2020}
Meiyi Ma, Ji~Gao, Lu~Feng, and John Stankovic.
\newblock Stlnet: Signal temporal logic enforced multivariate recurrent neural
  networks.
\newblock {\em Advances in Neural Information Processing Systems},
  33:14604--14614, 2020.

\bibitem[MSA16]{scvx}
Yuanqi Mao, Michael Szmuk, and Behçet Açıkmeşe.
\newblock Successive convexification of non-convex optimal control problems and
  its convergence properties.
\newblock In {\em 2016 IEEE 55th Conference on Decision and Control (CDC)},
  pages 3636--3641, 2016.

\bibitem[PDN21]{Puranic2021}
Aniruddh~G Puranic, Jyotirmoy~V Deshmukh, and Stefanos Nikolaidis.
\newblock Learning from demonstrations using signal temporal logic in
  stochastic and continuous domains.
\newblock {\em IEEE Robotics and Automation Letters}, 6(4):6250--6257, 2021.

\bibitem[Pnu77]{pnueli1977}
Amir Pnueli.
\newblock The temporal logic of programs.
\newblock {\em 18th Annual Symposium on Foundations of Computer Science (sfcs
  1977)}, pages 46--57, 1977.

\bibitem[RDM{\etalchar{+}}14]{Raman2014}
Vasumathi Raman, Alexandre Donzé, Mehdi Maasoumy, Richard~M. Murray, Alberto
  Sangiovanni-Vincentelli, and Sanjit~A. Seshia.
\newblock Model predictive control with signal temporal logic specifications.
\newblock In {\em 53rd IEEE Conference on Decision and Control}, pages 81--87,
  2014.

\bibitem[Ric05]{richards2005}
Arthur Richards.
\newblock Trajectory optimization using mixed-integer linear programming.
\newblock 05 2005.

\bibitem[VCDJS17]{vazquez2017}
Marcell Vazquez-Chanlatte, Jyotirmoy~V Deshmukh, Xiaoqing Jin, and Sanjit~A
  Seshia.
\newblock Logical clustering and learning for time-series data.
\newblock In {\em Computer Aided Verification: 29th International Conference,
  CAV 2017, Heidelberg, Germany, July 24-28, 2017, Proceedings, Part I 30},
  pages 305--325. Springer, 2017.

\bibitem[WBT21]{wang2021}
Changhao Wang, Jeffrey Bingham, and Masayoshi Tomizuka.
\newblock Trajectory splitting: A distributed formulation for collision
  avoiding trajectory optimization.
\newblock In {\em 2021 IEEE/RSJ International Conference on Intelligent Robots
  and Systems (IROS)}, pages 8113--8120, 2021.

\bibitem[YYW{\etalchar{+}}19]{yang2019}
Tao Yang, Xinlei Yi, Junfeng Wu, Ye~Yuan, Di~Wu, Ziyang Meng, Yiguang Hong,
  Hong Wang, Zongli Lin, and Karl~H. Johansson.
\newblock A survey of distributed optimization.
\newblock {\em Annual Reviews in Control}, 47:278--305, 2019.

\end{thebibliography}

\end{document}